\documentclass[11pt]{article}
\usepackage{amsfonts,amsmath,latexsym,verbatim,amscd}
\usepackage{amsthm,amsmath,amssymb,amscd}
\usepackage{amssymb}
\usepackage{graphics}

\def\a{{\alpha}}

\def\b{{\partial}}

\def \o{\omega}

\def\e{{\varepsilon}}

\def \ch{ {\cosh } \, \o}
\def \cdh{ {\cosh}^2 \, \o}

\def \sh{{\sinh} \, \o}

\newcommand{\R}{{\mathbb R}}

\newcommand{\N}{{\mathbb N}}

\newcommand{\Z}{{\mathbb Z}}

\newtheorem{lemma}{Lemma}[section]
\newtheorem{proposition}{Proposition}[section]
\newtheorem{theorem}{Theorem}[section]
\newtheorem{remark}{Remark}[section]

\newtheorem{definition}{Definition}[section]

\title{Higher genus Riemann minimal surfaces}

\author{Laurent Hauswirth \thanks{Email: Laurent.Hauswirth@univ-mlv.fr}\\ Universit\'e
de Marne la Vall\'ee  \and Frank
Pacard \thanks{Email: pacard@univ-paris12.fr, membre de l'Institut Universitaire de France}\\
Universit\'e Paris 12}

\date{October 07, 2005}

\begin{document}

\maketitle

\section{Introduction}

B. Riemann \cite{riemann} has constructed a one parameter family of
non congruent singly periodic minimal surfaces which are foliated by
circles (or straight lines). Each member of this family is a
periodic embedded minimal surface in $\R^3$ with infinitely many
parallel ends.

\medskip

Even though the classification of genus zero, embedded minimal
surfaces is not complete, W. H. Meeks J. Perez and A. Ros
\cite{MPR1}, \cite{MPR2}, \cite{MPR3} have made progress concerning
the question of the uniqueness of the Riemann examples in the class
of genus zero embedded minimal surfaces which  have an
infinite number of ends. They conjecture in \cite{MPR2} that every
embedded minimal surface of finite genus and with infinite number of
ends is asymptotic (away from a compact piece) to some "middle"
planar end and to two halves of Riemann example which are referred
to as "limit ends".

\medskip

In this paper we construct such surfaces. More precisely, we have
the~:
\begin{theorem}
Given $k =1, \ldots, 37$, there exists a one parameter family of
properly embedded minimal surfaces of genus $k$ with two limit ends  asymptotic
to half Riemann surfaces.
\end{theorem}

We briefly explain the idea behind the proof, this will give further
information about the surfaces constructed. In 1981, C. Costa
\cite{costa1}, \cite{costa2} found a genus one minimal properly
embedded surface, with three ends, two of which are asymptotic to
catenoidal ends and the third one being asymptotic to a plane.
Later, D. Hoffman and W. H. Meeks \cite{HM1}, \cite{HM2} have found
for every genus $k \geq 2$ a minimal surface with finite topology,
two catenoidal ends and one planar end.

\medskip

Minimal surfaces belonging to Riemann's family, once they are
normalized so that their planar ends are horizontal and at distance
$1$ one from each other, depend on a parameter (basically the value
of the horizontal flux). As this parameter tends to $0$, the members
of this family can be understood as infinitely many horizontal
planes linked by slightly bent catenoid.

\medskip

The main idea behind our construction is to replace one of these
"slightly bent" catenoid by one minimal surface which belongs to the
Costa-Hoffman-Meeks family of minimal surfaces. Our main result says
that this can be construction is successful provided one can bend
the upper and lower end of the genus $k$ Costa-Hoffman-Meeks
surface. Thanks to the moduli space theory for minimal surfaces with
catenoidal ends and a nondegeneray result by S. Nayatani
\cite{nayatani}, we are able to show that the bending of the ends of
the genus $k$ Costa-Hoffman-Meeks surface is possible for $1 \leq k
\leq 37$.

\medskip

The paper is organized as follows~: In Section 2, we give a
description of the Costa-Hoffman-Meeks  minimal surfaces and we
proceed with the deformation of the top and bottom ends of such
surfaces. In Section 3, we describe an isothermal parametrization of
Riemann surface, we also obtain some important expansions of pieces
of Riemann's surfaces as the flux becomes vertical. Section 4 is
devoted to the study of the mapping properties of the Jacobi
operator about a half Riemann surface as the flux becomes vertical.
In Section 5, we apply the implicit function theorem to perturb a
half Riemann surface, we obtain an infinite dimensional family of
minimal surfaces which are asymptotic to a half Riemann surface. In
Section 6, we perturb the Costa-Hoffman-Meeks surface using again
the implicit function theorem, we again obtain an infinite
dimensional family of minimal surfaces which have two boundaries and
one horizontal end. In the last section, we explain how the boundary
data of the minimal surfaces constructed in Section 4 and Section 5
an be chosen so that the union of theses forms a smooth minimal
surface with fixed genus and two limit ends.

\section{The Costa-Hoffman-Meeks' family of minimal surfaces}

C. Costa \cite{costa1}, \cite{costa2} and later on D. Hoffman and W.
H. Meeks \cite{HM1}, \cite{HM2} have described, for $k \geq 1$, a
properly embedded minimal surface of genus $k$ with three ends. More
precisely, for each $k \geq 1$, there exists $M_k$ a complete
properly embedded minimal surface of genus $k$ and three ends which,
after suitable rotation and translation, enjoys the following
properties~:

\begin{enumerate}
\item[(i)]
The surface $M_k$ has one planar end $E_m$ asymptotic to the $x_3=0$
plane, one top end $E_t$ asymptotic to the upper end of a catenoid
with $x_3$-axis of revolution and one bottom end $E_b$ asymptotic to
the lower end of a catenoid with $x_3$-axis of revolution. The
planar end $E_m$ is located in between the two catenoidal ends. \\

\item[(ii)]
The surface $M_k$ is invariant under the action of the rotation of
angle $\frac{2\pi}{k+1}$ about the $x_3$-axis, it is also invariant
under the action of the symmetry with respect to the $x_2 =0$ plane.
Finally, it is invariant under the action of the composition of a
rotation of angle $\frac{\pi}{k+1}$ about the $x_3$-axis and the
symmetry with respect to the $x_3 =0$ plane.\\

\item[(iii)]
The surface $M_k$ intersects the $x_3 =0$ plane in $k+1$ straight
lines, which intersect at equal angles $\frac{\pi}{k+1}$ at the
origin. The intersection of $M_k$ with the plane $x_3 = cte \, (\neq
0)$ is a single Jordan curve. The intersection of $M_k$ with the
upper half space $x_3 > 0$ (resp. with the lower half space $x_3
 < 0$) is topologically an open
annulus.
\end{enumerate}
The surface $M_k$ will be referred to as the "genus $k$
Costa-Hoffman-Meeks surface". Observe that, when $k$ is even the
surface $M_k$ is also invariant under the action of the rotation of
angle $\pi$ about the $x_2$-axis.

\medskip

The main purpose of this section is to explain how the genus $k$
Costa-Hoffman-Meeks surface $M_k$ can be deformed into a smooth one
parameter family of minimal surfaces $M_k ( \xi )$, for $\xi \in
(-\xi_0, \xi_0)$ and $\xi_0 > 0$ small enough, which are not
embedded anymore, are invariant under the action of the symmetry
with respect to the $x_2=0$ plane, have one horizontal end
asymptotic to the $x_3 = 0$ plane and have two catenoidal type ends
which are (up to some translations) respectively asymptotic to the
upper end and the lower end of a catenoid whose axis of revolution
is directed by $\sin \xi \, e_1  + \cos \xi \, e_3 $. The
construction of $M_k (\xi )$ will be a simple consequence of the
moduli space theory as described in \cite{perez-ros},
\cite{Kus-Maz-Pol} or \cite{Jle}. It also relies on a nondegeneracy
assumption which is known to be true when $k \leq 37$, thanks to
result of S. Nayatani \cite{nayatani}.

\medskip

Given $k \geq 1$, we start with a local description of the surface
$M_k$ near its ends and in particular we describe coordinates which
will be used to define some weighted spaces of functions on $M_k$.
The planar end $E_m$ of the surface $M_k$ can be parameterized by
\[
X_m (x) : =   \left( \frac{x}{|x|^2}, u_m (x) \right) \in {\mathbb
R}^3
\]
where $x  \in B_{r_0}(0) -\{0\} \subset {\mathbb R}^2$ and where the
function $u_m$ tends to $0$ as $x$ tends to $0$. This reflects the
fact that the middle end of $M_k$ is asymptotic to the horizontal
plane. Here $r_0 >0$ is fixed large enough.

\medskip

Recall that, for surfaces parameterized by
\[
x  \longrightarrow \left(
\frac{x}{|x|^2}, u (x) \right) \in {\mathbb R}^3
\]
the minimal surface equation reads
\begin{equation}
|x|^4 \, \mbox{div} \, \left( \frac{\nabla u}{(1+ |x|^4 \, |\nabla
u|^2)^{1/2}} \right) =0.
\label{eq:2.1}
\end{equation}
The function $u_m$ is (by definition) a solution of this equation
and it turns out that $u_m$, which is {\it a priori} only defined in
$B_{r_0} (0)- \{ 0 \}$, can be extended smoothly to $B_{r_0}$. We
shall make use of this fact, which follows from elliptic regularity
theory, without further comment. Observe that $u_m (x) = {\mathcal
O} (|x|)$ near $0$, however, given the symmetry with respect to the
rotation of vertical axis and angle $\frac{2\pi}{k+1}$, one checks
that $u_m (x) = {\mathcal O} (|x|^{k+1})$ near $0$. Indeed, since
$u_m$ solves (\ref{eq:2.1}), the leading term in the expansion of
$u_m$ in powers of $|x|$ is necessarily a harmonic function which is
invariant under the action of a rotation of angle $\frac{2\pi}{k+1}$
hence, in polar coordinates, it is a linear combination of the
functions $(r, \theta ) \longrightarrow r^{k+1}\, e^{\pm i \, (k+1)
\, \theta}$.

\medskip

We now turn to the description of the top end of $M_k$ (the
description of the bottom end will follow at once using the
invariance of the surface $M_k$ by the symmetries which are
described in (ii)). As already mentioned, the top end is asymptotic
to a catenoid with vertical axis of revolution. We use
\[
X_c (s, \theta) : = ( \cosh s \, \cos \theta, \cosh s \, \sin
\theta, s) \in {\mathbb R}^3
\]
as a parametrization of the (standard) catenoid $C$ with $x_3$-axis
of revolution. The unit normal vector field about $C$ is chosen to
be
\[
N_c (s, \theta) : = \frac{1}{\cosh s} \, ( \cos \theta,  \sin
\theta,  - \sinh s).
\]
Up to some dilation, we can assume that the top end $E_t$ of the
surface $M_k$ is asymptotic to some translated copy of the catenoid
parameterized by $X_c$ in the vertical direction. Therefore, $E_t$
can be parameterized by
\[
X_t : =  X_c  + w_t   \, N_c  + \sigma_t \, e_3
\]
for $(s, \theta) \in (s_0, \infty) \times S^1$, where the function
$w_t$ tends to $0$ as $s$ tends to $\infty$ and $\sigma_t \in
{\mathbb R}$. Again, $w_t$ tends to $0$ as $s$ tends to $\infty$,
reflecting the fact that the end $E_t$ is asymptotic to the standard
catenoid translated by $\sigma_t \, e_3$.

\medskip

We recall that the surface parameterized by $X : =  X_c  + w \, N_c$
is minimal if and only if the function $w$ satisfies the minimal
surface equation which,  for normal graphs over a catenoid, can be
expanded in powers of $w$ (and its partial derivatives) as
\begin{equation}
\frac{1}{\cosh^2 s} \, \left( \left( \partial_s^2 +
\partial_\theta^2 + \frac{2}{\cosh^2 s} \right) \, w +  Q_2
\left( \frac{w}{\cosh s}\right) + \cosh s \, Q_3 \left(
\frac{w}{\cosh s}\right)\right) =0 \label{eq:2.2}
\end{equation}
Here $Q_2$ and $Q_3$ are nonlinear second order differential
operators which  satisfy
\begin{equation}
\| Q_j (v_2) - Q_j (v_1) \|_{{\mathcal C}^{0, \alpha} ((s,s+1)
\times S^1)} \leq c \, \left( \sup_{i=1,2} \| v_i \|_{{\mathcal
C}^{2, \alpha} ((s,s+1) \times S^1)} \right)^{j-1}  \,  \| v_2 -v_1
\|_{{\mathcal C}^{2, \alpha} ((s,s+1) \times S^1)} \label{eq:ZZZ}
\end{equation}
for all $s \in {\mathbb R}$ and all $v_1, v_2$ such that $\| v_i
\|_{{\mathcal C}^{2, \alpha} ((s,s+1) \times S^1)} \leq 1$. The
important fact is that the constant $c>0$ does not depend on $s$.
The proof of this expansion can be easily adapted from the proof of
the corresponding expansion for higher dimensional catenoids which
is provided in \cite{Fak-Pac}, a complete (short) proof is given in
the Appendix A.

\medskip

The function $w_t$ is (by definition) a solution of (\ref{eq:2.2}).
Given the symmetry with respect to the rotation of vertical axis and
angle $\frac{2\pi}{k+1}$, one checks that $w_t$ is in fact bounded
by a constant times $e^{-(k+1) s}$. Indeed, just observe that, in
the expansion of $w_t$ in powers of $e^{-s}$, the leading term is
harmonic (on the cylinder ${\mathbb R} \times S^1$) and invariant
under the action of the rotation on $S^1$ by the angle
$\frac{2\pi}{k+1}$, hence it has to be a linear combination of the
functions $(s, \theta )  \longrightarrow e^{-(k+1) s} \, e^{\pm i \,
(k+1) \, \theta}$.

\medskip

Similarly, we define $X_b$ to parameterize the lower end $E_b$ of
the surface $M_k$ so that
\[
X_b : =  X_c  - w_b   \, N_c  - \sigma_b \, e_3
\]
for $(s, \theta) \in (-  \infty, - s_0,) \times S^1$, where the
function $w_b$ tends to $0$ as $s$ tends to $-\infty$ and $\sigma_b
\in {\mathbb R}$. Again, $w_b$ tends to $0$ as $s$ tends to
$-\infty$, reflecting the fact that the end $E_b$ is asymptotic to
the standard catenoid translated by $\sigma_b \, e_3$. Granted the
symmetries of the surface $M_k$, there is an obvious relation
between $X_t$ and $X_b$. Indeed, starting from the parametrization
of $E_t$ which we compose by a rotation of angle $\frac{\pi}{k+1}$
about the $x_3$-axis and a symmetry with respect to the $x_3=0$ one
finds a parametrization of $E_b$. This implies that  $\sigma_b =
\sigma_t$ and also that
\[
w_b ( s, \theta ) =  - w_t \left( -s, \theta - \frac{\pi}{k+1}
\right).
\]

For all $r < r_0$ and $s > s_0$, we define
\begin{equation}
M_k(s,r) : = M_k  - \left( X_t ((s ,\infty) \times S^1) \cup X_b
((-\infty , -s) \times S^1) \cup X_m (B_{r }(0)) \right)
\label{eq:2.3}
\end{equation}

The parametrizations of the three ends of $M_k$ induce a
decomposition of $M_k$ into slightly overlapping components as
follows : A compact piece $M_k (s_0+1, r_0/2)$ and three noncompact
pieces $X_t ((s_0 ,\infty) \times S^1)$, $X_b ((-\infty, -s_0)
\times S^1)$ and $X_m (B_{r_0} (0))$. We are now in a position to
define the weighted spaces of functions on $M_k$.
\begin{definition}
Given $\ell \in {\mathbb N}$, $\alpha \in (0,1)$ and $\delta , \nu
\in {\mathbb R}$, the space ${\mathcal C}^{\ell , \alpha}_{\delta ,
\nu} (M_k)$ is defined to be the space of functions in ${\mathcal
C}^{\ell , \alpha}_{loc} (M_k)$ for which the following norm is
finite
\[
\begin{array}{rllll}
\| w \|_{{\mathcal C}^{\ell , \alpha}_{\delta , \nu} (M_k)}&  : = &
\| w \|_{{\mathcal C}^{\ell , \alpha} (M_k (s_0+1, r_0/2))} + \| \,
|\cdot|^{-\nu} \, w \circ X_m \|_{{\mathcal C}^{\ell , \alpha}  (B_{r_0}(0))} \\[3mm]
& + & \sup_{s \geq s_0} e^{-\delta s} \,  \left( \| w \circ X_t
\|_{{\mathcal C}^{\ell , \alpha} ((s,s+1) \times S^1)} + \| w \circ
X_b \|_{{\mathcal C}^{\ell , \alpha} ((-s-1,-s) \times S^1)} \right)
\end{array}
\]
and which are invariant under the action of the symmetry with
respect to the $x_2=0$ plane, i.e.  $w (p)  = w(\bar p)$  for all
$p\in M_k$, where $\bar p : = (x_1, -x_2, x_3)$ if $p = (x_1, x_2,
x_3) $. \label{de:2.1}
\end{definition}

The Jacobi operator about $M_k$ is defined by
\[
{\mathbb L}_{M_k} : = \Delta _{M_k} + | A_{M_k}|^2
\]
where $|A_{M_k}|$ is the norm of the second fundamental form on
$M_k$. Granted the above defined spaces, one can check that~:
\[
\begin{array}{rccccllll}
L_\delta : &  {\mathcal C}^{2, \alpha}_{\delta ,0} (M_k) &
\longrightarrow & {\mathcal C}^{0, \alpha}_{\delta - 2 , 4}
(M_k)\\[3mm]
& w &  \longmapsto & {\mathbb L}_{M_k} \, (w)
\end{array}
\]
is a bounded linear operator. The subscript $\delta$ is meant to
keep track of the weighted space over which the Jacobi operator is
acting. Observe that, in the weights of the target space, there is a
loss of $2$ in the weight parameter at the ends $E_t$ and $E_b$, and
there is a gain of $4$ in the weight parameter at the end $E_m$.
This follows at once from the expression of the Jacobi operator at
the ends in the above defined coordinates. Alternatively, this can
also be seen by linearizing the nonlinear equation (\ref{eq:2.1}) at
$u = 0$ which provides the expression of the Jacobi operator about
the plane
\[
{\mathbb L}_{{\mathbb R}^2} : = |x|^4 \, \Delta
\]
and by linearizing the nonlinear equation (\ref{eq:2.2}) at $w = 0$
which provides the expression of the Jacobi operator about the
standard catenoid
\[
{\mathbb L}_C : = \frac{1}{\cosh^2 s} \, \left( \partial_s^2+
\partial_\theta^2 + \frac{2}{\cosh^2 s}\right)
\]
Since the Jacobi operator about $M_k$ is asymptotic to ${\mathbb
L}_{{\mathbb R}^2}$ at $E_m$ and is asymptotic to ${\mathbb L}_C$ at
$E_t$ and $E_b$, this explains the loss of $2$ in the weight
parameter $\delta$ and the gain of $4$ in the weight parameter
$\nu$.

\medskip

This being understood, we now recall the notion of nondegeneracy
\cite{Kus-Maz-Pol} which is classically used in this context~:
\begin{definition}
The surface $M_k$ is said to be nondegenerate if $L_\delta$ is
injective for all $\delta < -1$. \label{de:2.2}
\end{definition}

The mapping properties of the operator $L_\delta$ depends crucially
on the choice of $\delta$. It follows from the general theory of
such operators that $L_\delta$ has closed range and is Fredholm
provided $\delta \notin {\mathbb Z}$. Moreover, a duality argument
(in weighted Lebesgue spaces !) implies that
\[
\left( L_\delta \quad \mbox{is injective} \right)  \quad
\Leftrightarrow \quad  \left( L_{-\delta} \quad \mbox{is
surjective}\right)
\]
provided $\delta \notin {\mathbb Z}$. This kind of analysis is by
now standard and has been applied to variety of problems.  We refer
to \cite{Mel} for references to the general theory and we refer to
\cite{Jle} for references to the theory in the specific context of
minimal hypersurfaces with catenoidal type ends. Also, we have the~:
\begin{proposition}
Assume that $M_k$ is nondegenerate and $\delta \in (1, 2)$. Then the
operator $L_\delta$ is surjective. Moreover the kernel of $L_\delta$
is $4$-dimensional.
\end{proposition}
One has to keep in mind that, in the definition of the weighted
spaces, we have imposed the invariance under some symmetry and that,
in addition, we have implicitly asked that the middle end of the
surface remain asymptotic to a horizontal plane. This explains why
the dimension of the kernel is only equal to $4$ and not equal to $9
(= 3 \times $ the number of ends) as is usually the case when no
symmetries are imposed.

\medskip

Recall that a smooth one parameter group of isometries containing
the identity generates a Jacobi field i.e. a solution of the
homogeneous problem ${\mathbb L}_{M_k} \, w =0$. We now define $4$
of these Jacobi fields and we also provide there expansion at the
ends of $M_k$. Let $N$ denote a unit normal vector field on $M_k$
(for example, we agree that the orientation is chosen so that $N
\sim e_3$ at $E_m$). We will denote by
\[
\Phi^{0,+}(p) : =  N (p) \cdot e_3,
\]
the Jacobi field generated by the one parameter group of vertical
translations. Observe that
\begin{equation}
\begin{array}{rlrllll}
\Phi^{0,+} & = & - \tanh s  + {\mathcal O}((\cosh s)^{-k-2}) & \quad
\mbox{at} \qquad E_t \\[3mm]
\Phi^{0,+} & = &  \tanh s + {\mathcal O}((\cosh s)^{-k-2})  & \quad
\mbox{at} \qquad  E_b
\end{array}
\label{eq:2.4}
\end{equation}
while $ \Phi^{0,+} = 1 + {\mathcal O} (|x|^{2k+4})$ at $E_m$. We
will denote by
\[
\Phi^{0,-} (p)  := N(p) \cdot p
\]
the Jacobi field generated by the one parameter group of dilations.
Observe that
\begin{equation}
\begin{array}{rlllll}
\Phi^{0,-}  & = & 1 - s  \, \tanh s  + {\mathcal O}((\cosh
s)^{-k-1}) &
\quad \mbox{at} \qquad E_t \\[3mm]
\Phi^{0,-} & =&  s \, \tanh s  - 1 + {\mathcal O}((\cosh s)^{-k-1})
& \quad \mbox{at} \qquad E_b
\end{array}
\label{eq:2.5}
\end{equation} while $ \Phi^{0,-} = {\mathcal O}
(|x|^{k+1})$ at $E_m$. We denote by
\[
\Phi^{1,+} (p) : = N(p) \cdot e_1
\]
the Jacobi field generated by the one parameter group of
translations along the $x_1$-axis. Observe that
\begin{equation}
\begin{array}{rlrllll}
\Phi^{1,+} & = & \frac{1}{\cosh s} \, \cos \theta + {\mathcal
O}((\cosh s)^{-k-2}) & \quad \mbox{at} \qquad E_t \\[3mm]
\Phi^{1,+} & = & - \frac{1}{\cosh s} \, \cos \theta + {\mathcal
O}((\cosh s)^{-k-2}) & \quad \mbox{at} \qquad E_b
\end{array}
\label{eq:2.6}
\end{equation} while $\Phi^{0,-} = {\mathcal O}
(|x|^{k+2}) $ at $E_m$. Finally, we denote by
\[
\Phi^{1,-} (p)  := N(p) \cdot (e_2 \times p)
\]
the Jacobi field generated by the one parameter group of rotation
about the $x_2$-axis. Observe that
\begin{equation}
\begin{array}{rlrllll}
\Phi^{1,-} & = &  (\frac{s}{\cosh s} + \sinh s ) \, \cos \theta +
{\mathcal O}((\cosh s)^{-k-1}) & \quad \mbox{at} \qquad E_t \\[3mm]
\Phi^{1,-} & = & - (\frac{s}{\cosh s} + \sinh s )  \, \cos \theta +
{\mathcal O}((\cosh s)^{-k-1}) & \quad \mbox{at} \qquad E_b
\end{array} \label{eq:2.7}
\end{equation} while $\Phi^{1,-} = \frac{x_1}{|x|^2} +
{\mathcal O} (|x|^{2k+3})$ at $E_m$.

\medskip

Observe that all these globally defined Jacobi fields are invariant
under the action of the symmetry with respect to the $x_2=0$ plane
and that there are in addition three other Jacobi fields which are
not invariant under this symmetry, namely the Jacobi field
associated to the group of translation along the $x_2$-axis and the
Jacobi field corresponding to the one parameter group of rotation
about the $x_1$-axis and the Jacobi field corresponding to the one
parameter group of rotation about the $x_3$-axis.

\medskip

With these notations, we define the deficiency space
\[
{\mathcal D} : =  \mbox{Span} \{ \chi_t \, \Phi^{j,\pm }, \chi_b \,
\Phi^{j, \pm} \qquad  : \qquad  j=0, 1\}
\]
where $\chi_t$ is a cutoff function  which is identically equal to
$1$ on $X_t ((s_0+1,\infty)\times S^1)$, identically equal to $0$ on
$M_k -X_t ((s_0,\infty)\times S^1)$ and which satisfies $\chi_t (p)
= \chi_t(\bar p)$ (so that it is invariant under the action of the
symmetry with respect to the $x_2 =0$ plane).  We also define
$\chi_b (\cdot) : = \chi_t (- \, \cdot)$. Clearly
\[
\begin{array}{rcccllll}
\tilde L_\delta :  & {\mathcal C}^{2, \alpha}_{\delta ,0} (M_k)
\oplus {\mathcal D} & \longrightarrow & {\mathcal C}^{0,
\alpha}_{\delta-2 , 4} (M_k) \\[3mm]
& w &  \longmapsto & {\mathbb L}_{M_k} \, (w)
\end{array}
\]
is a bounded linear operator.

\medskip

The linear decomposition Lemma proved in \cite{Kus-Maz-Pol} for
constant mean curvature surfaces or in \cite{Jle} for minimal
hypersurfaces can be adapted to our situation and we get the~:
\begin{proposition}
Assume that $M_k$ is nondegenerate and that $\delta \in (-2, -1)$.
Then the operator $\tilde L_\delta$ is surjective and has a kernel
of dimension $4$. \label{pr:2.2}
\end{proposition}

We are interested in ${\mathcal M}$ the space of all minimal
surfaces (not necessarily embedded) which are close to $M_k$, have
$2$ catenoidal ends, one horizontal planar end and which are
invariant under the action of the symmetry with respect to the
$x_2=0$ plane. The moduli space theory developed in
\cite{Kus-Maz-Pol} for constant mean curvature surfaces or in
\cite{Jle} for minimal hypersurfaces can be adapted to our framework
and as a corollary of Proposition~\ref{pr:2.2}, we conclude that,
close to $M_k$, the space ${\mathcal M}$ is a smooth manifold of
dimension $4$, provided $M_k$ is nondegenerate. Moreover, the
elements of the kernel of $\tilde L_\delta$ span the tangent space
to ${\mathcal M}$. Therefore, in order to understand the space
${\mathcal M}$ in a neighborhood of $M_k$, we just need to
understand the elements which span the kernel of $\tilde L_\delta$
since this will provide the set of parameters which are needed to
describe ${\mathcal M}$ in a neighborhood of $M_k$.

\medskip

It should be clear that the functions $\Phi^{0,\pm}$ and
$\Phi^{1,+}$ belong to ${\mathcal C}^{2, \alpha}_{\delta ,0} (M_k)
\oplus {\mathcal D}$ and hence we already have $3$ linearly
independent elements of the kernel of $\tilde L_\delta$. Observe
that $\Phi^{1,-}$ fails to belong to the kernel of $\tilde L_\delta$
since it is not bounded at $E_m$ (and in fact blows up like
$|x|^{-1}$ as $x$ tends to $0$). Thus, we are left to understand the
behavior of a nonzero element $\Phi \in {\mathcal C}^{2,
\alpha}_{\delta ,0} (M_k) \oplus {\mathcal D}$ which belongs to the
kernel of $\tilde L_\delta$ but does not belong to $\mbox{Span}
\{\Phi^{0,\pm} , \Phi^{1,+}\}$, and hence $\Phi  \neq 0$. Without
loss of generality (i.e. taking suitable linear combination of
$\Phi$ with $\Phi^{0,\pm}$ and $\Phi^{1,+}$) we can assume that the
expansion of $\Phi$ at $E_t$ is given by
\[
\Phi =  a_t \, \Phi^{1,-} + {\mathcal O} ((\cosh s)^\delta)
\]
and that the expansion of $\Phi$ at $E_b$ is given by
\[
\Phi =  a_b \, \Phi^{1,-} + b_b \, \Phi^{1, +} + c_b \, \Phi^{0,+} +
d_b \, \Phi^{0,-} +  {\mathcal O} ( (\cosh s)^\delta )
\]

Given a function $\Psi$ defined on $M_k$, we set
\[
W (\Psi) : = \lim_{s\rightarrow \infty} \, \lim_{r\rightarrow 0}
\int_{M_k(s,r)} \, (\Phi \, {\mathbb L}_{M_k} \, \Psi - \Psi \,
{\mathbb L}_{M_k} \, \Phi) \, dvol_{M_k}
\]
where we recall that $M_k(s ,r )$ has been defined in
(\ref{eq:2.3}).

\medskip

Since ${\mathbb L}_{M_k}\Phi =0$ and ${\mathbb L}_{M_k} \, \Phi^{0, +}
=0$, we can use the divergence theorem together with the expansions
(\ref{eq:2.4})-(\ref{eq:2.5}) to get
\[
0 = W(\Phi^{0,+}) = 2 \, \pi \, d_b.
\]
Similarly, using the fact that ${\mathbb L}_{M_k} \, \Phi^{0, -} =0$
together with the expansions (\ref{eq:2.4})-(\ref{eq:2.5}), we get
\[
0 = W (\Phi^{0,-}) = - 2 \, \pi \, c_b
\]
Next, using the fact that ${\mathbb L}_{M_k} \, \Phi^{1, -} =0$
together with the expansions (\ref{eq:2.6})-(\ref{eq:2.7}), one
finds that
\[
0 = W (\Phi^{1,-}) =  - \, \pi \, b_b
\]
Finally, using the fact that ${\mathbb L}_{M_k} \, \Phi^{1, +} =0$
together with the expansions (\ref{eq:2.6})-(\ref{eq:2.7}),  we have
\[
0 = W (\Phi^{1,-}) =  \pi \, (a_b  - a_t).
\]
Therefore, we conclude that $b_b = c_b = d_b=0$ and also that
$a_b=a_t$. Now, if we had $a_t = 0$, then we would also have $a_b=0$
and hence we would conclude that $\Phi \in {\mathcal C}^{2,
\alpha}_\delta (M_k)$. But, in this case nondegeneracy implies that
$\Phi = 0$, which is clearly a contradiction since we have assumed
that $\Phi\neq 0$. Therefore, we conclude that $a_t \neq 0$. In
other words, there exists an element of the kernel of $\tilde
L_\delta$ which at $E_t$ (and in fact also at $E_b$) is asymptotic
to the Jacobi field associated to the rotation of the catenoidal
ends of $M_k$, leaving the middle end horizontal.

\medskip

Applying (an elaborate version of) the implicit function theorem as
in \cite{Kus-Maz-Pol} and \cite{Jle}, we see that this Jacobi field
is integrable. This shows that there exists in ${\mathcal M}$ a one
parameter family of minimal hypersurfaces $(M_k (\xi))_\xi$, for
$\xi$ close to $0$, such that $M_k(0)=M_k$ and the catenoidal upper
end of $M_k (\xi)$ is asymptotic to the end of a catenoid whose axis
of revolution is directed by $\sin \xi \, e_1 + \cos \xi \, e_3$.
Observe that, the surface $M_k (\xi)$ is well defined up to a
translation in the $x_2=0$ plane and up to a dilation. In
particular, we can require that the upper end of $M_k(\xi)$ is
asymptotic to a translated and rotated version of the (standard)
catenoid and also require that the middle end $E_m(\xi)$ is
asymptotic to the $x_3=0$ plane.

\medskip

If $R_\xi$ denotes the rotation of angle $\xi$ about the $x_2$-axis,
the upper end $E_t(\xi)$ of $M_k(\xi)$ can be parameterized by
\begin{equation}
X_{t,\xi} = R_\xi  \, ( X_c + w_{t, \xi} \, N_c) + \sigma_{t ,\xi}
\, e_3  + \varsigma_{t ,\xi} \, e_1 \label{eq:paramt}
\end{equation}
where the function $w_{t, \xi}$ and $\sigma_{t ,\xi}, \varsigma_{t
,\xi} \in {\mathbb R}$ depend smoothly on $\xi$ and satisfy  $w_{t,
0} = w_{t}$, $\sigma_{t ,0} =\sigma_t$ and $\varsigma_{t ,0} = 0$.
More precisely, if follows from the application of the implicit
function theorem that
\[
|\sigma_{t ,\xi } - \sigma_{t}| + |\varsigma_{t ,\xi }| + \| w_{t,
\xi} - w_{t}\|_{{\mathcal C}^{2, \alpha}_{-2} ((s_0,+\infty) \times
S^1)} \leq \, c \, |\xi|
\]
Application of the flux formula \cite{KKS} shows that the lower end
of $M_k(\xi)$ is, up to a translation, asymptotic to the lower end
of the same (standard) catenoid. In particular, the lower end $E_b
(\xi)$ of $M_k(\xi)$ can be parameterized by
\begin{equation}
X_{b,\xi} = R_\xi  \, ( X_c - w_{b, \xi} \, N_c) - \sigma_{b ,\xi}
\, e_3 - \varsigma_{b ,\xi} \, e_1 \label{eq:paramb}
\end{equation}
where the function $w_{b, \xi}$ and  $\sigma_{b  ,\xi}, \varsigma_{b
, \xi} \in {\mathbb R}$ depend smoothly on $\xi$ and satisfy  $w_{b,
0} = w_{b}$, $\sigma_{b ,0} = \sigma_b$ and $\varsigma_{b ,0} = 0$.
Again, we also have
\begin{equation}
|\sigma_{b ,\xi } - \sigma_b | + |\varsigma_{b ,\xi}|  + \| w_{b,
\xi} - w_{b}\|_{{\mathcal C}^{2, \alpha}_{-2} ((-\infty , - s_0)
\times S^1)} \leq \, c \, |\xi| \label{eq:11211}
\end{equation}
Finally, up to now the surface $M_k(\xi)$ is defined up to a
translation along the $x_1$-axis but we can help eliminate this
confusion by requiring that
\[
\varsigma_{b,\xi} = \varsigma_{t,\xi}
\]

To summarize, we have obtained the~:
\begin{theorem}
Assume that $M_k$ is nondegenerate. Then, there exists $\xi _0 >0$
and a smooth one parameter family of minimal hypersurfaces $(M_k
(\xi))_\xi$ in ${\mathcal M}$, for $\xi \in (- \xi_0, \xi_0)$, such
that $M_k(0)=M_k$ and the upper (resp. lower) catenoidal end of $M_k
(\xi)$ is, up to a translation along its axis, asymptotic to the
upper (resp. lower) end of the standard catenoid whose axis of
revolution is directed by $\sin \xi \, e_1 + \cos \xi \, e_3$.
\end{theorem}

Observe that, when $k$ is even, the surface $M_k$ is symmetric with
respect to the rotation of angle $\pi$ about the $x_2$-axis and one
can prove, even if we will not need this information, that the
surfaces $M_k (\xi)$ can be defined to enjoy the same symmetry.

\medskip

On each $M_k(\xi)$ one can define weighted spaces as in
Definition~\ref{de:2.1} and also define the corresponding notion of
nondegeneracy. The Jacobi operator about $M_k (\xi)$ will be denoted
by ${\mathbb L}_{M_k(\xi)}$. We define
\[
\begin{array}{rccccllll}
L_{\xi , \delta} :  & {\mathcal C}^{2, \alpha}_{\delta ,0}
(M_k(\xi))& \longrightarrow & {\mathcal C}^{0, \alpha}_{\delta-2 ,
4} (M_k(\xi))\\[3mm]
& w & \longmapsto & {\mathbb L}_{M_k(\xi)} \, (w)
\end{array}
\]
It is easy to check that, reducing $\xi_0$ if this is necessary, all
the surfaces $M_k (\xi)$ are nondegenerate and hence we have the~:
\begin{proposition}
Assume that $M_k$ is nondegenerate and choose $\delta \in (1, 2)$.
Then (reducing $\xi_0$ if this is necessary) the operator $L_{\xi,
\delta}$ is surjective and has a kernel of dimension $4$. Moreover,
there exists $G_{\xi , \delta }$ a right inverse for $L_{\xi ,
\delta}$ which depends smoothly on $\xi$ and in particular whose
norm is bounded uniformly as $|\xi| < \xi_0$. \label{pr:lk}
\end{proposition}
For example, a right inverse $G_{\xi,\delta}$ which depends smoothly
on $\xi$ can be obtained by a simple perturbation argument starting
from a right inverse $G_{\delta, 0}$ and reducing $\xi_0$ if this is
necessary.

\medskip

The purpose of the next Lemma is to write a portion of the upper and
lower ends of the surface $M_k (\xi)$ as vertical graphs over the
horizontal plane. It is clear that the ends $E_t$ and $E_b$ of $M_k$
can be written, at least away from a compact set, as vertical graphs
over the horizontal plane $x_3=0$. This is not true anymore for the
ends $E_{t} (\xi)$ and $E_{b} (\xi)$ of $M_k (\xi)$, when $\xi \neq
0$ but this property will remain true for the piece of $M_k (\xi)$
we are interested in, namely the piece of $M_k (\xi)$ which
corresponds to $s \sim \pm \frac{1}{2} \, \log \e$ in the
parametrization given in (\ref{eq:paramt}) and (\ref{eq:paramb}). We
have the~:

\begin{lemma}
\label{dl3} There exists $\e_0 >0$ such that, for all $\e \in (0,
\e_0)$ and all $|\xi |\leq \e$ an annular part of $E_t(\xi)$ and
$E_b(\xi)$ in $M_k (\xi)$ can be written as vertical graphs over the
horizontal plane for the functions
\[
\begin{array}{rlr}
\bar U^0_{t} (r,\theta ) & = & \sigma_{t, \xi} + \ln (2 r) + \xi \,
r \, \cos \theta + {\mathcal O} (\e)
\\[3mm]
\bar U_{b}^0 (r , \theta ) & = & - \sigma_{b, \xi} - \ln (2r)  + \xi
\, r \, \cos \theta + {\mathcal O} (\e)
\end{array}
\]
Here $(r, \theta)$ are polar coordinates in the $x_3=0$ plane. The
functions ${\mathcal O}(\e)$ are defined in the annulus $B_{4 \,
\e^{-1/2}} - B_{\e^{-1/2} / 4}$ and are bounded in ${\mathcal
C}^\infty_b$ topology by a constant (independent of $\e$) times
$\e$, where partial derivatives are computed with respect to the
vector fields $r \, \partial_r$ and $\partial_\theta$.
\end{lemma}
{\bf Proof~:} Elementary computations writing
\[
X_{t , \xi} (s, \tilde \theta) = (r  \, \cos \theta, r \, \sin
\theta , U_{t}^0 (r, \theta))
\]
when $s \sim - \frac{1}{2} \, \log \e$. \hfill $\Box$

\medskip

We end this section by recalling the result of S. Nayatani
\cite{nayatani} which states that $M_k$ is nondegenerate for all $k
\leq 37$. More precisely we have the~:
\begin{theorem}
Assume that $k \leq 37$, then any bounded Jacobi field on $M_k$ is a
linear combination of  $N\cdot e_j$ and $N \cdot (p\times e_3)$.
\end{theorem}
In order to apply this result, just observe that, according to this
result, when $k \leq 37$, the only Jacobi field which decays at all
ends is $ N\cdot (p\times e_3)$. However, this Jacobi field is not
invariant with respect to the action of $p \rightarrow \bar p$,
hence $M_k$ is nondegenerate in the sense defined in
Definition~\ref{de:2.2}.

\section{Riemann minimal surface}

B.  Riemann \cite{riemann} has discovered a one parameter family of
periodic minimal surfaces embedded in ${\R}^3$ which are foliated by
circles (and straight lines). Each element of this family has
infinitely many planar ends, is topologically a cylinder ${\mathbb
R} \times S^1 $ and in fact is conformal to the cylinder ${\mathbb
R}  \times S^1$ with infinitely many points  $(p_i)_{i \in \Z}$
removed in a periodic way, each of these points corresponds to one
of the planar ends of the surface.

\medskip

Recall that (up to some  dilation and some rigid motion) we can
parameterize a fundamental piece of Riemann's surface by
\begin{equation}
(t, \theta) \longrightarrow (c(t) + R(t) \, \cos \theta, R(t)  \,
\sin \theta, t) \label{eq:3.1} \end{equation} where $t \in (- t_\e,
t_\e)$, $\theta \in S^1$ and where the functions $c$ and $R$ are
determined by
\begin{equation}
(\partial_t R)^2 +1 = R^2 + \e^2 \, R^4
\label{eq:3.11}
\end{equation}
and
\[
\partial_t c = \e \, R^2
\]
Here $\e >0$ is a parameter. We shall normalize the solutions of
these ordinary differential equations by asking that $R(0)
>0$, $\partial_t R (0) =0$ and $c(0)=0$, and naturally, $R$ is a
nonconstant smooth solution of (\ref{eq:3.11}). Even though $R$ and
$c$ both depend on $\e$, we shall not make this dependence explicit
in the notation. It is easy to check that the functions $R$ and $c$
blow up in finite time $t_\e <\infty$ and that
\[
\ell_\e = \lim_{t\rightarrow t_\e} ( c(t)-R(t))
\]
exists. Riemann surface $R_\e$ is then obtained by translation of
the fundamental piece by $2 \, (\ell_\e \, e_1 + t_\e \, e_3) \,
\mathbb Z$.

\medskip

A conformal parametrization of Riemann surfaces had already been
considered by M. Shiffman \cite{shiffman} and has been generalized
by L. Hauswirth \cite{hauswirth}. Granted the above parametrization,
in order to define this conformal parametrization, it is enough to
look for a function $(t,y) \longrightarrow \psi (t,y)$ such that
\begin{equation}
X_\e (t, y) : = (c(t) + R(t) \, \cos (\psi (t,y)), R(t) \, \sin
(\psi (t,y)) , t) \label{eq:Xe} \end{equation}
is a conformal
parametrization. This leads to the first order differential system
\[
\partial_t \psi  = \e \, R \, \sin \psi
\]
which come from the equation $\partial_t X_\e \cdot \partial_y X_\e
=0$ and
\[
(\partial_y \psi )^2 = 1 + \e^2 \, R^2 \, (1 + \cos^2 \psi ) + 2 \,
\e \, \partial_t R \, \cos \psi
\]
which comes from the requirement that $|\partial_t X_\e |^2 =
|\partial_y X_\e |^2$. One checks easily (from a direct computation)
that the integrability condition $\partial_y \, ( \partial_t \psi) =
\partial_t \, (\partial_y \psi)$ is fulfilled.

\medskip

We define the real valued function $\o$ by
\[
\partial_t X_\e = (\sinh \o \, \cos \psi , \sinh \o \, \sin \psi, 1)
\qquad \mbox{and} \qquad \partial_y X_\e = ( - \cosh \o \, \sin \psi
, \cosh \o \, \cos \psi, 0)
\]
In particular
\begin{equation}
\cosh \o  =  R \, \partial_y \psi \label{eq:wRpsi}
\end{equation}

With these notations, the first fundamental form about the surface
parameterized by $X_\e$ reads
\[
ds^2 = \cosh^2 \o \, (dt \otimes dt + dy \otimes dy)
\]
and, if we define the normal vector field by
\begin{equation}
N_\e : = \frac{1}{\cosh \o} \, \left(\cos \psi, \sin \psi , - \sinh
\o \right),
\label{eq:Ne}
\end{equation}
the second fundamental form about the surface parameterized by
$X_\e$ is then given by
\[
h =  \partial_t \o \, dt \otimes dt - \partial_y \psi \, dy \otimes
dy - 2 \, \partial_t \psi \, dt \otimes dy
\]
Observe that $\partial_y (\partial_t  X) = \partial_t (\partial_y
X)$ and hence
\begin{equation}
\partial_y \o = - \partial_t \psi \qquad \mbox{and} \qquad
\partial_t \o = \partial_y \psi \label{eq:wpsi}
\end{equation}

It will be convenient to define
\begin{equation}
a : = \frac{\partial_t \psi}{\cosh \o}  \qquad \mbox{and} \qquad b =
\frac{\partial_y \psi}{\cosh \o} \label{eq:ab}
\end{equation}
With these notations, the Jacobi operator about Riemann's surface
$R_\e$ is given by
\[
{\mathbb L}_{R_\e} =  \frac{1}{\cosh^2 \omega} \, \left(
\partial_t^2 +
\partial_y^2 + 2 \, (a^2+b^2) \right)
\]

Observe that it follows from (\ref{eq:wRpsi}) that $b = \frac{1}{R}$
and hence the equation satisfied by $b$ reads
\[
(\partial_t b)^2+ b^4 =  \e^2 + b^2
\]
Moreover, $\partial_y b=0$ since $R$, and hence $b$, does not depend
on $y$.

\medskip

It should be clear that (\ref{eq:ab}) together with (\ref{eq:wpsi})
yields $\partial_t a + \partial_y b =0$ and hence the function $a$
does not depend on $t$. It remains to find the ordinary differential
equation satisfied by $a$. We have
\[
\partial_y (\ch \, a) - \partial_t (\ch
\, b)=0 .
\]
Hence
\begin{equation}
(a^2+ b^2) \, \sh = \partial_y a - \partial_t b \label{eq:1.7}
\end{equation}
Taking the derivative of (\ref{eq:1.7}) with respect to $y$ and
using the fact that the function $b$ does not depend on $t$, we get
\begin{equation}
\partial_y^2 a = - a \, (a^2+ b^2) +\frac{a}{a^2+b^2} \,
( (\partial_y a)^2 - (\partial_t b)^2) \label{eq:1.8}
\end{equation}
In other words
\[
\partial_y \, \left( \frac{ (\partial_y a)^2 -
(\partial_t b)^2}{a^2+b^2} + a^2 \right) =0
\]
Hence, we conclude that the function
\[
(t,y) \longrightarrow \frac{ (\partial_y a)^2 - (\partial_t
b)^2}{a^2+b^2} + a^2
\]
only depends on $t$. Taking this information into account in
(\ref{eq:1.8}) we conclude that
\[
\partial_y^2 a + 2 \, a^3 + \alpha (t) \, a = 0
\]
where $\alpha$ is a function which only depends on $t$, and hence
$\alpha$ has to be constant.  Therefore
\[
(\partial_y a )^2 +  a^4 + \alpha \, a^2 - \beta =0
\]
for some fixed constant $\beta$. Inserting these into
(\ref{eq:1.8}), we conclude that
\begin{equation}
(\partial_y a )^2 +  a^4 +  a^2 - \e^2 =0 \label{eq:1.11}
\end{equation}
The function $y \longrightarrow a(y)$ and $t \longrightarrow b(t)$
are defined up to some translation in the $y$ or $t$ variables. In
particular, we can require that  $a$ (resp. $b$) takes its maximal
value at $y=0$ (resp. $t=0$). Finally, observe that the functions
$a$ is periodic, we will denote by $y_{\e}$ its least period.
Finally, we extend the function $b$ to be a $2\, t_\e$-periodic
function.

\medskip

It is also easy to check that, as $\e$ tends to $0$ the functions
$\e^{-1} \, a$, $b$ and their derivatives remain uniformly bounded.
Indeed, we have on the one hand
\begin{equation}
(\partial_y a )^2 \leq \e^2 \qquad \mbox{and} \qquad 2\, a^2 \leq
\sqrt{1+4\e^2}-1 \label{eq:GF}
\end{equation}
and on the other hand
\begin{equation}
(\partial_t b )^2 \leq \e^2 + \frac{1}{4} \qquad \mbox{and} \qquad 2
b^2 \leq 1 + \sqrt{1+4\e^2}. \label{eq:equb}
\end{equation}
A simple application of Ascoli-Arzela's Theorem implies the~:
\begin{lemma}
As $\e$ tends to $0$, the sequence of functions $(b)_{\e >0}$
converges uniformly on compacts to the function
\[
t \longrightarrow \frac{1}{\cosh t}
\]
and the sequence of functions $(\e^{-1} \, a)_{\e >0}$ converges
uniformly to the function
\[
y \longrightarrow \cos y
\]
\end{lemma}

We can also obtain the expansion of the period $2 \, t_\e$ of the
function $b$. Indeed, we have the formula
\[
t_\e = \int_{0}^{\zeta_\e} \frac{d\zeta}{\sqrt{1 + \zeta^2- \e^2 \,
\zeta^4}}
\]
where $0 < \zeta_\e$ is the largest root of $\e^2 \, \zeta^4 =
\zeta^2 +1$.  It is easy to check that
\[
t_\e = -  \log \e + {\mathcal O} (1)
\]
as $\e$ tends to $0$.

\medskip

We claim that
\[
\frac{1}{\sqrt{1+ 4 \e^2}} \leq  \left(
\frac{y_\e}{2\pi} \right)^2 \leq 1
\]
Indeed, write $a(y) = a(0)  \, \cos v$ where $2 \, a(0)^2 +1  =
\sqrt{1+ 4 \e^2}$ and using (\ref{eq:GF}) we get
\[
(\partial_y v)^2 = 1 + a(0)^2 \, ( 1+ \cos^2 v)
\]
from which it follows that $1 \leq \partial_y v \leq 1+ 2 \,
a(0)^2$. Integration of this inequalities from $0$ to $y_\e$ yields
the required inequalities since $v(y_\e) =2 \, \pi$.

\medskip

In the next Lemma, we give precise expansions of the functions $R$
and $c$ when $t \in ( - t_\e, t_\e )$.
\begin{lemma}
\label{dl} For $\e >0$ small enough, we have
\begin{equation}
R (t)= \cosh t + {\mathcal O}(\e^2 \cosh^3 t) \label{eq:estR}
\end{equation}
and
\[
c(t) = \e  \,\left( \frac{t}{2} + \frac{\sinh (2t)}{4} \right) +
{\mathcal O}(\e^3 \cosh^4 t)
\]
when $t \in [ - t_\e + 1 , t_\e -1]$.
\end{lemma}
{\bf Proof~:} We define the function $t \longrightarrow v(t)$ such
that $R (t)= R(0) \, \cosh v(t)$ and $v(0)=0$. It follows from
(\ref{eq:3.11}) that
\[
(\partial_t v)^2 = 1 + \e^2  \, R^2 (0) \, (1 + \cosh ^2 v)
\]

\noindent

Now, as long as  $t \leq v(t) \leq t + c$ (where $c>0$ is some fixed
constant), we can estimate $(\partial_t v)^2 = 1 + {\mathcal O} (\e
^2 \cosh^2 t)$ and hence we conclude that $v (t)=  t + {\mathcal O}
(\e ^2 \cosh^2 t)$. We remark {\it a posteriori} that $t \leq v(t)
\leq t+ c$ holds for $t \in [0, - \ln \e -1]$ provided $c >0$ is
chosen large enough. The first estimate then follows at once. The
second estimate follows directly from
\[
\partial_t a =  \e \, R^2
\]
once the first estimate has been established. \hfill $\Box$

\begin{remark}
\label{re:remm} As a Corollary of the proof of this Lemma, observe
that
\[
 b (t)  \, \cosh t \,\leq \frac{1}{R(0)}
\]
for all $t \in [-t_\e, t_\e]$ since we always have $v (t)\geq t$.
Since $R(0)$ converges to $1$ as $\e$ tends to $0$ and using the
fact that $b$ is even, this yields a uniform upper bound for $b$ as
$\e$ tends to $0$.
\end{remark}

The purpose of the next Lemma is to write the pieces of $R_\e$ at
height $t \sim - \frac{1}{2} \log \e$  (resp. at height $t \sim
\frac{1}{2} \log \e$) as a vertical graph over the horizontal plane
for some function $t_t$ (resp. $t_b$). But before doing so we first
dilate the surface $R_\e$ by some factor $(1+ \gamma)$ and we next
translate this dilated surface along the $x_1$-axis by $\varsigma$,
so that the fundamental piece of this surface is now parameterized
by
\[
(t, y) \longrightarrow ( \varsigma + (1+\gamma) \, (c(t) + R(t) \,
\cos \psi(t,y)), (1+ \gamma) \, R(t) \, \sin \psi(t,y) , (1+\gamma)
\, t)
\]
We consider the change of coordinates~:
\[
( r \cos \theta , r \, \sin \theta  ) = ( \varsigma + (1+ \gamma) \,
( c(t) +R(t) \, \cos \psi (t,y) ), (1+\gamma) \, R (t) \sin \psi
(t,y))
\]
where as before, $(r, \theta)$ are polar coordinates in the $x_3=0$
plane. Obviously this change of coordinates is not valid everywhere
but we are only interested in the range $t \sim \pm \frac{1}{2} \,
\log \e$ where the change of coordinates holds.
\begin{lemma}
\label{graphe} Assume that $|\varsigma |\leq 1$ and also assume that
$|\gamma |\leq \frac{1}{2}$. Then the following expansion holds
\begin{equation}
U_t^0 (r,\theta) = (1+\gamma) \, \log \left(\frac{2
r}{1+\gamma}\right) - \frac{\e}{2} \, r \, \cos \theta - (1+\gamma)
\, \frac{\varsigma}{r} \, \cos \theta + {\mathcal O} (\e)
\label{eq:ttop}
\end{equation}
and
\begin{equation}
U_b^0 (r,\theta) = -  (1+\gamma) \, \log \left(\frac{2
r}{1+\gamma}\right)  - \frac{\e}{2} \, r \, \cos \theta + (1+\gamma)
\, \frac{\varsigma}{r} \, \cos \theta + {\mathcal O} (\e)
\label{eq:tbot}
\end{equation}
when $ r \sim \e^{-1/2}$. Here the functions ${\mathcal O} (\e)$ are
smooth functions which are defined in the annulus $B_{4 \,
\e^{-1/2}} - B_{\e^{-1/2}/4}$ and are bounded by a constant
(independent of $\e$) times $\e$ in ${\mathcal C}^\infty_b$
topology, where partial derivatives are understood with respect to
the vector fields $r \, \partial_r $ and $\partial_\theta$. In
addition all these estimates hold uniformly in $\sigma$ and
$\gamma$, provided $|\varsigma|\leq 1$ and $|\gamma | \leq
\frac{1}{2}$. \label{le:t1t2}
\end{lemma}

\section{The Jacobi operator about Riemann's surface}

We keep the notations of the previous section. Recall that the
Jacobi operator about Riemann's surface is given by
\begin{equation}
{\mathbb L}_{R_\e} : =  \frac{1}{\cosh^2 \omega} \, \left(
\partial^2_t +\partial_{y}^2 + 2 \, (a^2 + b^2) \right)
\label{eq:4.0}
\end{equation}
Obviously the mapping properties of this operator translate into
the mapping properties of the operator
\begin{equation}
L_\e = \partial^2_t + \partial_{y}^2  + 2 \, (a^2 + b^2)
\label{eq:2.233}
\end{equation}

We define, for all $\e \geq 0$ the operator
\[
D_\e : = \partial_{y}^2  + 2 \, a^2
\]
which acts on functions of $y$ which are $y_\e$ periodic and even.
This operator is clearly elliptic and self adjoint and hence has
discret spectrum $(\lambda_{i})_{i\geq 0}$. Since we only consider
even functions, each eigenvalue is simple and we can arrange the
eigenvalues so that $\lambda_{i} < \lambda_{i+1}$. The corresponding
eigenfunctions are denoted by $f_{i}$ and are normalized so that
\[
\int_0^{y_\e} f_{i}^2 \, dy =1
\]
Even though we have not make this explicit in the notations, the
eigendata of $D_\e$ do depend on $\e$ since $a$ does. It is easy to
check that, as $\e$ tends to $0$, the $\lambda _{i}$ converge to
$i^2$. We will not need this result but rather the simpler~:
\begin{lemma}
\label{lemme7} The following estimate holds
\begin{equation}
\lambda _{i} \geq  i^2 + 1 -  \sqrt{1+4\e^2} . \label{eq:4.1}
\end{equation}
\end{lemma}
{\bf Proof :} The assertion follows from the variational
characterization of the eigenvalues
\[
\lambda_{i} = \sup_{\mbox{codim} \, E = i} \, \inf_{f \in E, \, \| f
\|_{L^2}=1} \int_0^{y_\e} \left( (\partial_{y} f)^2 - 2 \, a^2 \,
f^2 \right) \, dy
\]
together with the fact that $ y_\e \leq 2 \pi$ and $0 \leq 2 a^2
\leq \sqrt{1+4\e^2} -1$. \hfill $\Box$

\medskip

For each $\e >0$, the family $\{ f_{i} \}_{i \in \N }$ is a Hilbert
basis of the space $H$ of $L^2$-integrable functions which are even
and $y_\e$-periodic. We consider the eigenfunction decomposition of
a function $(t,y) \longrightarrow v(t,y)$, which is $y_\e$-periodic
and even in the $y$ variable,
\[
v (t,y) = \sum ^{\infty}_{i=0} v_i(t) \, f_{i}(y).
\]
This decomposition induces a decomposition of the operator $L_\e$
into the sequence of ordinary differential operators
\[
L_{\e,i} = \partial_{y}^2 + 2 \, b ^2  - \lambda _{i}.
\]
It follows from the result of Lemma~\ref{lemme7} and from the
estimate (\ref{eq:4.1}) that
\begin{equation}
2 \, b^2  - \lambda _{i} \leq 2 \, \sqrt{1 + 4 \, \e^2} - i^2
\end{equation}

Let $S^1 (\tau)$ denote the circle of radius $\tau$. The previous
estimate immediately implies the following injectivity result, by
maximum prinnciple~:
\begin{lemma}
\label{lemme8} Assume that $\e \in (0, \sqrt{\frac{3}{4}})$ and $t_0
< t_1$. Let $v$ be a solution of
\[
L_{\e} \, v =0
\]
which is defined on $[t_0, t_1] \times S^1 (\tau_\e)$ and satisfies
$v(t_0, \cdot ) = v(t_1, \cdot )=0$. Further assume that
\[
\int_{0}^{y_\e} v(t,y) \, f_{i} (y)\, dy =0
\]
for all $t \in [t_0, t_1]$ and $i=0,1$. Then $v=0$.
\end{lemma}
The constant $\sqrt{\frac{3}{4}}$ is not optimal but this will be
sufficient for our purposes since our aim is to use the result for
small values of $\e$. We now define weighted H\"older spaces which
will turn to be useful for the understanding the mapping properties
of the operator $L_\e$ as the parameter $\e$ tends to $0$.

\begin{definition}
Given $\ell \in \N $, $\a \in (0,1)$, $\mu \in \R$ and a closed
interval $I \subset {\R}$, we define the space ${\mathcal C}^{\ell,
\a}_{\mu} ( I \times S^1(\tau)) $ to be the space of functions $u
\in {\mathcal C}^{\ell, \a}_{loc} ( I \times S^1(\tau) )$ for which
the following norm
\[
\| u \|_{{\mathcal C}^{\ell, \a}_{\mu}} =  \| e^{-\mu t} \, u
\|_{{\mathcal C}^{\ell, \a}(I \times S^1(\tau))},
\]
is finite.
\end{definition}

We set
\[
\tau_\e = \frac{y_\e}{2\pi}
\]
It should be obvious that
\[
\begin{array}{rlllll}
{\mathcal C}^{2, \a}_{\mu}([t_0, \infty) \times S^1(\tau_\e)) &
\longrightarrow & {\mathcal C}^{0, \a}_{\mu} ([t_0, \infty) \times
S^1(\tau_\e) )\\[3mm]
w & \longmapsto & L_\e \, w
\end{array}
\]
for any $\mu \in {\R}$ and $t_0 \in {\R}$. We prove that, provided
the parameter $\mu$ is suitably chosen, there exists a right inverse
for $L_\e$ whose norm is uniformly bounded as $\e$ tends to $0$ and
independently of $t_0 \in {\R}$. This is the content of the
following~:
\begin{proposition}
\label{inverse} Fix $\mu \in (-2, -1)$. Then, there exists $\e_0 >0$
and, for all $\e \in (0, \e_0)$, for all $t_0 \in {\R}$, there
exists an operator
\[
G_{\e,t_0} : {\mathcal C}^{0, \a}_{\mu}([t_0, \infty) \times S ^1
(\tau_\e)) \longrightarrow  {\mathcal C}^{2, \a}_{\mu} ([t_0,
\infty) \times S ^1 (\tau_\e)),
\]
such that for all $g \in  {\mathcal C}^{0, \a}_{\mu}([t_0, \infty)
\times S^1 (\tau_\e) )$, the function $v :=G_{\e, t_0}(g)$ solves
\[
\left\{ \begin{array}{rlll}
L_\e \, v &  =  & g  \qquad \mbox{in} \quad  [t_0, \infty) \times S^1 ( \tau_\e)  \\[3mm]
        v & \in & \mbox{Span} \, \{f_{0}, f_{1}\}  \quad \mbox{on} \quad  \{t_0\}
\times S^1 ( \tau_\e)
\end{array}\right.
\]
Moreover,
\[
||G_{\e, t_0}(g)||_{ {\mathcal C}^{2, \a}_{\mu} } \leq c \,
\|g\|_{{\mathcal C}^{0, \a}_{\mu}},
\]
for some constant  $c>0$  which is independent of $\e  \in (0,
\e_0)$ and also independent of $t_0 \in {\mathbb R}$.
\end{proposition}
{\bf Proof~:} We decompose $f$ into
\[
g = g_0 \, f_{0} + g_1 \, f_{1} + \bar g
\]
where $\bar g (t, \cdot )$ is $L^2$ orthogonal to $f_{0}$ and
$f_{1}$ for each $t$. For the sake of implicity in the notations, we
shall not mention the parameter $\tau_\e$  and write $S^1$ instead
of $S^1(\tau_\e)$. Observe that, as $\e$ tends to $0$, $\tau_\e$
tends to $1$.

\medskip

{\bf Step 1.} We show that, for each $t_1 > t_0 +1$ it is possible
to solve
\[
L_{\e} \, \bar v = \bar g ,
\]
on $S^1 \times [t_0,t_1]$ with $\bar v ( t_0 , \cdot ) =\bar v (t_1,
 \cdot ) = 0$. This just follows from the result of Lemma~\ref{lemme8}
which states that, restricted to the set of functions $L^2$
orthogonal to $f_{0}$ and $f_{1}$ for each $t$, the operator $L_\e$
is injective.

\medskip

We claim that, provided $\e$ is chosen small enough, there exists
a constant $c >0$ such that
\[
\sup_{[t_0, t_1] \times S^1} e^{-\mu t} \, |\bar v|\leq c \,
\sup_{[t_0, t_1] \times S^1} e^{-\mu t} \, |\bar g|
\]
The proof of this fact is by contradiction. If this were false,
there would exist a sequence $(\e_n)_n$ tending to $0$, sequences
$(t_{0,n})_n$ and $(t_{1,n})_n$ such that $t_{0,n} \leq t_{1,n}$, a
sequence of functions $(\bar g_n)_n$ and a sequence of solutions
$(\bar v_n)_n$ such that
\[
\sup_{[t_{0,n}, t_{1,n}] \times S^1 } e^{-\mu t} \, |\bar g_n| =1
\]
and
\[
\lim_{n\rightarrow +\infty} A_n : = \sup_{[t_{0,n}, t_{1,n}] \times
S^1} e^{-\mu t} \, |\bar v_n| =  \infty
\]
We denote by $(t_n,y_n) \in [t_{0,n}, t_{1,n}] \times S^1$ a point
where $A_n$ is achieved. We define the function $\tilde v_n$ by
\[
\tilde v_n (t,y)= \frac{1}{A_n} \, e^{- \mu t_n}\, \bar v_n (t_n +
t, y).
\]

Observe that elliptic estimates imply that
\begin{equation}
\sup e^{-\mu t} |\nabla \bar v_n| \leq c \, (1+ A_n)
\end{equation}
and, since $\bar v_n$ vanishes on the boundaries of $ [t_{0,n},
t_{1,n}] \times S^1$, this in turn implies that the sequences $(t_n
- t_{0,n})_n$ and $(t_{1,n} - t_n)_n$ remain bounded away from $0$.

\medskip

Without loss of generality, we can assume that the sequence $(t_n -
t_{0,n})_n$ (resp. ($t_{1,n} - t_0)_n$) converges to $\bar t_0 \in
([-\infty, 0)$ (resp. to $\bar t_1 \in (0, +\infty]$). We denote by
$I = (\bar t_0, \bar t_1)$.

\medskip

Up to a subsequence, we can assume, without loss of generality that
the sequence of functions $(\tilde v_n)_n$ converges on compacts to
a nontrivial function $\tilde v$ defined on $I \times S^1$. This
follows at once from Ascoli-Arzela theorem, once it is observed that
the sequence of functions $(\tilde v_n)_n$ is uniformly bounded (by
$t \longrightarrow e^{\mu t}$) and, by elliptic regularity theory,
the sequence of functions $(\nabla \tilde v_n)_n$ is also uniformly
bounded (by $t \longrightarrow c \, e^{\mu t}$). We now derive some
properties of the limit function $\tilde v$. These properties are
all inherited by $\tilde v$ from similar properties which hold for
the functions $\bar v_n$.

\medskip

First, $\tilde v(t, \cdot )$ is $L^2$ orthogonal to the constant
function and the function $y \rightarrow \cos y$ for each $t \in I$.
Next, $\tilde v$ is equal to $0$ on $\{\tilde t_0\} \times S^1$ and
on $\{\tilde t_1\} \times S^1$ if either $\tilde t_0
> - \infty$ or $\tilde t_1 < +\infty$. Also
\begin{equation}
\sup_{I \times S^1} e^{-\mu t} \, |\tilde v| = 1 \label{eq:4.5}
\end{equation}
Finally, $\tilde v$ is either a solution of
\begin{equation}
\left( \partial^2_{t} + \partial_{y}^2 + \frac{2}{\cosh^2 (\cdot +
\tilde t)} \right) \, \tilde v = 0 \label{eq:4.6}
\end{equation}
for some $\tilde y \in {\R}$ or is a solution of
\[
(\partial^2_{t}  + \partial_{y}^2 ) \, \tilde v = 0.
\]

To reach a contradiction we consider the eigenfunction
decomposition of $\tilde v$
\[
\tilde v (t,y) = \sum_{j=2}^\infty a_j(t) \, \cos (j\, y) .
\]

When $\bar t_0 = -\infty$, observe that the function $a_j$ is either
blowing up like $t \longrightarrow e^{- j t}$ or decaying like $y
\longrightarrow e^{j t}$. The choice of $\mu \in (-2, -1)$ implies
that $a_j$ decays exponentially at $-\infty$. Multiplying the
equation (\ref{eq:4.6}) by $a_j \, e_j$ and integrating by parts
over $I$ (all integrations are justified because $a_j$ decays
exponentially at both $\pm \infty$ if either $\bar t_0 = -\infty$ or
$\bar t_1 = +\infty$), we get either
\[
\int_{-\infty}^{+\infty} ( |\partial_t a_j|^2 + j^2 \,a_j^2 ) \, dt
= \int_{-\infty}^{+\infty} \frac{2}{\cosh^2 (t + \tilde t)} \, a_j^2
\, dt
\]
or
\[
\int_{-\infty}^{+\infty} ( |\partial_t a_j|^2 + j^2 \, a_j^2 ) \, dt
= 0
\]
In either case, we obtain $a_j=0$ which clearly contradicts
(\ref{eq:4.5}).

\medskip

Since we have reached a contradiction, the proof of the claim is
complete. Once the claim is proven, we can use once more elliptic
estimates and Ascoli-Arzela theorem to pass to the limit as $t_1$
tends to $+\infty$ in a sequence of solutions which are defined on
$[t_0, t_1] \times S^1$. This proves the existence of a solution of
\[ L_\e \, \bar v  = \bar f \] which is defined in $S^1 \times [t_0,
+ \infty)$ and which satisfies $\bar v (t_0, \cdot ) =0$. In
addition, we know that
\[
\sup_{[t_0, +\infty) \times S^1 } e^{-\mu t} \, |\bar v|\leq c \,
\sup_{[t_0, +\infty) \times S^1} e^{-\mu t} \, |\bar f|
\]
Using a last time elliptic estimates, we complete the proof of the
result in the case where the eigenfunction decomposition of $f$
does not involve $e_{\e,0}$ or $e_{\e,1}$.

\medskip

\noindent {\bf  Step 2.} Now we consider the case where the function
$g$ is collinear to $e_{0}$ and $e_{1}$, namely
\[
g(t,y)= g_0(t)\, f_{0}(y) + g_1(t)\, f_{1}(y)
\]
We extend the function $g$ to be equal to $0$ when $t \leq t_0$,
keeping the same notation. Given $t_1 > t_0$, we consider the
equation
\[
( \partial_{t}^2 + 2 \, b^2 - \lambda_{j} ) \, v_j = g_j
\]
in $(-\infty, t_1)$ with boundary data $v_j(t_1)= \partial_t v_j
(t_1)=0$. The existence of $v_j$ is standard. We claim that
\[
\sup_{(-\infty, t_1)} e^{- \mu t} \, |v_j| \leq c\, \sup_{\R} e^{-
\mu t} \, |g_j|
\]
for some constant which does not depend on $t_1$, provided $\e$ is
chosen small enough. As before, we argue by contradiction. Assume
that the claim is not true, there would exist a sequence $(\e_n)_n$
tending to $0$, sequences $(t_{0,n})_n$ and $(t_{1,n})_n$ such that
$y_{0,n} \leq t_{1,n}$, a sequence of functions $(g_{j,n})_n$ and a
sequence of solutions $(v_{j,n})_n$ such that
\[
\sup_{ (-\infty , t_{1,n}] \times S^1 } e^{-\mu t} \, |g_{j,n}| =1
\]
and
\[
\lim_{n\rightarrow +\infty} A_n : = \sup_{(-\infty , y_{1,n}] \times
S^1} e^{-\mu t} \, |v_{j,n}| = \infty
\]
We denote by $(t_n,y_n) \in (-\infty , t_{1,n}] \times S^1$ a point
where $A_n$ is achieved. Observe that, the solution $v_{j,n}$ is a
linear combination of the two solutions of the homogeneous problem
$L_{\e, j} \, w =0$ and these are known to be at most linearly
growing thanks to Jacobi fields coming from isometries. Hence the
above supremum is achieved. We define the function $\tilde v_{j,n}$
by
\[
\tilde v_{j,n} (t,y)= \frac{1}{A_n} \, e^{- \mu t_n}\, v_{j,n} (t_n
+ t , y).
\]
As above, one shows that the sequence $(t_{1,n} - t_n)_n$ remains
bounded away from $0$.

\medskip

Without loss of generality, we can assume that the sequence
($t_{1,n} - t_0)_n$ converges to  $\bar t_1 \in (0, +\infty]$. We
denote by $I = (-\infty, \bar t_1)$.

\medskip

As in Step 1,  we can also assume, without loss of generality that
the sequence of functions $(\tilde v_{j,n})_n$ converges on compacts
to a nontrivial function $\tilde v_j$ defined on $I \times S^1 $. We
now derive some properties of the limit function $\tilde v_j$.

\medskip

First, $\tilde v_j (t, \cdot)$ is $\tilde v$ is equal to $0$ on
$\{\tilde t_1\} \times S^1$ if   $\tilde t_1 < +\infty$. Also
\begin{equation}
\sup_{I \times S^1} e^{-\mu t} \, |\tilde v_j| = 1 \label{eq:4.9}
\end{equation}
Finally, $\tilde v_j$ is either a solution of
\begin{equation}
\left( \partial_t^2  - j^2 + \frac{2}{\cosh^2 (\cdot + \tilde t)}
\right) \, \tilde v_j = 0 \label{eq:4.10}
\end{equation}
for some $\tilde t \in {\R}$ or is a solution of
\begin{equation}
(\partial_t^2 - j^2 ) \, \tilde v_j = 0. \label{eq:4.11}
\end{equation}

When $n=0$, the solutions of (\ref{eq:4.10}) are linear combinations
of the functions
\[
t \longrightarrow \tanh ( t + \tilde t) \qquad \mbox{and} \qquad t
\longrightarrow (t + \tilde t) \, \tanh (t + \tilde t) -1
\]
and when $n=1$ they are linear combinations of the following
functions
\[
t \longrightarrow \frac{1}{\cosh (t + \tilde t)}  \qquad \mbox{and}
\qquad t \longrightarrow \frac{(t + \tilde t)}{\cosh (t + \tilde t)}
+ \sinh (t+\tilde t)
\]

Finally, when $n=0$ the solutions of (\ref{eq:4.11}) are linear
combinations of the functions
\[
t \longrightarrow 1 \qquad \mbox{and} \qquad t \longrightarrow t
\]
and when $n=1$ they are linear combinations of the functions
\[
t \longrightarrow e^t \qquad \mbox{and} \qquad t \longrightarrow
e^{-t}
\]
To reach a contradiction we observe that all the solutions of
these two equations, when $j=0,1$ are explicitly known and that
none of them satisfies (\ref{eq:4.9}) since we have chosen $\mu \in
(-2, -1)$.

\medskip

Since we have reached a contradiction, the proof of the claim is
complete. Once the claim is proven, we pass to the limit as $t_1$
tends to $+\infty$ in a sequence of solutions which are defined on
$(-\infty , t_1] \times S^1$. This proves the existence of a
solution of
\[
L_{\e,j} \, v_j  = g_j
\]
which is defined in $[t_0, +\infty) \times S^1$. In addition, we
know that
\[
\sup_{ [t_0, +\infty) \times S^1} e^{-\mu t} \, |v_j|\leq c \,
\sup_{[t_0, +\infty) \times S^1 } e^{-\mu t} \, |g_j|
\]
The estimate for the derivatives follow from standard elliptic
estimates.  \hfill $\Box$

\medskip

The following result is standard and left to the reader (a proof can
be found in \cite{Fak-Pac}).
\begin{lemma}
\label{poissonbis} There exists an operator
\[
{\mathcal P} : {\mathcal C}^{2, \a}(S ^1) \longrightarrow {\mathcal
C}^{2, \a}_{-2} ([0,+ \infty) \times S^1),
\]
such that for all $\varphi \in  {\mathcal C}^{2, \a}(S ^1 )$, with
$\varphi$ orthogonal to $1$ and $\theta \longrightarrow \cos \theta
$ in the $L^2$-sense and is an even function of $y \in S^1 $ the
function $w = {\mathcal P}(\varphi)$ solves
\[
\left\{ \begin{array}{rllllll} (\b ^{2}_{t} + \partial_\theta^{2} )
\, w
& =& 0 & \mbox{ in }& [0,+ \infty)  \times S^1  \\[3mm]
w & = & \varphi & \mbox{ on }& \{0 \} \times S^1
\end{array}\right.
\]
Moreover,
\[
\|{\mathcal P}(\varphi) \|_{{\mathcal C}^{ 2, \a}_{ -2} } \leq c  \,
\| \varphi \|_{{\mathcal C}^{2, \a}},
\]
for some constant  $c>0$.
\end{lemma}

\section{An infinite dimensional family of minimal surfaces which
are close to a half Riemann surface}

In this section we are interested in minimal surfaces which are
close to a half of Riemann's surface and have prescribed boundary.
We consider surfaces which are normal graphs over Riemann's surface.
More precisely, we consider the surface parameterized by
\[
Z_{\e, u} : = X_\e + u \, N_\e.
\]
where $X_\e$ and $N_\e$ have been defined in (\ref{eq:Xe}) and
(\ref{eq:Ne}) and, in the following Proposition, we give an
expansion of the mean curvature operator for this surface in terms
of the function $u$ and its partial derivatives.
\begin{proposition}
\label{quadratique} The surface parameterized by $Z_{\e, u}:= X_\e
+u \, N_\e$ is minimal if and only if the function $u$ is a solution
of
\[
L_\e \, u = (\cosh \o )^2 \, Q_\e \left( \frac{u}{\cosh \o } ,
\frac{\nabla u}{\cosh \o }, \frac{\nabla^2 u}{\cosh \o }\right)
\]
\noindent where $L_\e$ is the operator which has already been
defined in (\ref{eq:2.233}) and the nonlinear operator $Q_\e$
satisfies
\[
| Q_\e (v_2) - Q_\e (v_1)|_{{\mathcal C}^{0, \alpha} ((t,t+1) \times
S^1(\tau_\e))} \leq c \, \sup_{i=1,2} | v_i|_{{\mathcal C}^{2,
\alpha} ((t,t+1) \times S^1(\tau_\e))} \,  | v_2 -v_1|_{{\mathcal
C}^{2, \alpha} ((t,t+1) \times S^1(\tau_\e))}
\]
for all $v_1, v_2$ such that $| v_i|_{{\mathcal C}^{2, \alpha}
((t,t+1) \times S^1(\tau_\e))} \leq 1$. Here the constant $c >0$
does not depend on $t \in {\mathbb R}$, nor on $\e \in (0,1)$.
\end{proposition}
{\bf Proof :} We omit the indices $\e$ and $u$ for the sake of
simplicity in the notations. We have
\[
\partial_t  Z = \partial_t  X + \partial_t  u \, N +  u \, \partial_t  N
\qquad \mbox{and} \qquad \partial_y  Z = \partial_y  X + \partial_y
u \, N +  u \, \partial_y  N
\]
A simple computation shows that the coefficients of $g_u$, the first
fundamental form of the surface parameterized by $Z_u$, are given by
\[
\begin{array}{rllllll}
|\b_t Z|^2 & = &  \cdh + 2 \, b \, \ch \, u +
(\partial_t u)^2 + (a^2 + b^2) \, u^2 \\[3mm]
 \b_t Z \cdot \b_y Z & = &  2 \, a \, \ch \, u + \partial_t u
\, \partial_y u \\[3mm]
|\b_y Z|^2 & =&  \cdh - 2 \, b \, \ch \, u+ (\partial_y u) ^2 + (a^2
+ b^2) \, u^2
\end{array}
\]
Collecting these, we have the expansion of the determinant of $g_u$
\[
\begin{array}{rllll} |g_u| & = &  \displaystyle \cosh^4 \o \, \left(1 +
\frac{1}{\cosh^2 \o} \, \left( |\partial_t u|^2 + |\partial_y u|^2 -
2 \, (a^2+ b^2) \, u^2\right)  \right. \\[3mm]
& & \displaystyle \qquad \qquad  \quad \left. + P_3 \left(
\frac{u}{\ch }, \frac{ \nabla u}{\ch } \right) + P_4 \left(
\frac{u}{\ch }, \frac{ \nabla u}{\ch } \right) \right)
\end{array}
\]
where $P_i$ has coefficients which are bounded independently of $\e
\in (0,1)$ are is homogeneous of degree $i$. Here we have implicitly
used the fact that the functions $a$ and $b$ are uniformly bounded
when $\e \in (0,1)$.

\medskip

We consider the area energy
\[
A (u) := \int \sqrt{|g_u|} \, dt\, dy
\]
and the surface parameterized by $Z$ will be minimal if and only if
the first variation if $0$. This can be written as
\[
2 \, D A_{|u} (v) = \int \frac{1}{\sqrt{|g_u|}} D_u |g_u| \, (v) \,
dt \, dy
\]
Observe that
\begin{equation}
\begin{array}{rllll} \frac{1}{\sqrt{|g_u|}} \, D |g_u| _{|u}\, (v)
& = &  \displaystyle 2 \, ( \partial_t u \, \partial_t v
+ \partial_y u \, \partial_y v - 2 \, (a^2+ b^2) \, u \, v )  \\[3mm]
& & \displaystyle +  \ch  \, Q \left( \frac{u}{\ch }, \frac{ \nabla
u}{\ch } \right) \, v + \ch \, \tilde Q \left( \frac{u}{\ch },
\frac{ \nabla u}{\ch } \right) \, \partial_t v \\[3mm]
&  & + \displaystyle \ch  \, \hat Q \left( \frac{u}{\ch }, \frac{
\nabla u}{\ch } \right) \, \partial_y v
\end{array}
\label{eq:ac27}
\end{equation}
where the operator $Q$, $\tilde Q$ and $\hat Q$ enjoy properties
similar to the one enjoyed by $Q_\e$ in the statement of the result.

\medskip

The result then follows at once provided one notices that
\[
\left| \partial_t \cosh \o  \right| +\left|
\partial_y \cosh \o \right| \leq (|\o_t|+|\o_y|) \, \cosh \o  \leq  (|a|+ |b|)  \, \cosh ^2 \o \leq
c \, \cosh^2 \o,
\]
for some constant $c >0$ which does not depend on $\e \in (0,1)$.
This explain the $\cosh^2 \o$ in front of the nonlinearity $Q_\e$
whereas (\ref{eq:ac27}) would only suggest a $\cosh \o$. \hfill
$\Box$

\medskip

We consider the surface parameterized by $X_{\e}$ which we first
dilate by a factor $(1+\gamma)$ and then translate by
$\frac{\varsigma}{1+\gamma}$ along the $x_1$ axis and by $(1+\gamma)
\, \log (1+\gamma) + \sigma$ along the $x_3$-axis. This surface,
which will be referred to as $R_{\e} (\gamma, \sigma, \varsigma)$ is
parameterized by
\begin{equation}
(1+\gamma ) \, X_{\e} + \frac{\varsigma}{1+\gamma} \, e_1 +
((1+\gamma) \, \log (1+\gamma) + \sigma)\, e_3. \label{eq:cbos}
\end{equation}
The parameters $\gamma $ and $\varsigma$ are now chosen to satisfy
\[
| \gamma | +  \e^{1/2} \, |\varsigma| \leq \kappa \, \e
\]
for some constant $\kappa >0$ which will be fixed later on.

\medskip

Using the result of Lemma~\ref{le:t1t2}, we see that part of this
surface (basically the one at height $t \sim - \frac{1+\gamma}{2} \,
\log \e$ is a graph over the annulus  $B_{4 \, \e^{-1/2}} -
B_{\e^{-1/2}/4}$ in the $x_3=0$ plane for the function
\begin{equation}
(r, \theta )  \longrightarrow (1+\gamma) \, \log (2 r ) + \sigma -
\frac{\e}{2} \, r \, \cos \theta - \frac{\varsigma}{r} \, \cos
\theta + {\mathcal O} (\e) \label{eq:vgtt}
\end{equation}
which has been expanded in (\ref{eq:ttop}). We now truncate the
surface $R_{\e} (\gamma, \sigma , \varsigma)$ at the graph of the
curve $r = \frac{1}{2} \, \e^{-1/2}$ by the function defined in
(\ref{eq:vgtt}) and consider only the upper half  of this surface
which we will refer to as $R_{\e }^t (\gamma, \sigma , \varsigma )$.

\medskip

We are interested in normal graphs over the surface $R_{\e }^t
(\gamma, \sigma , \varsigma)$ which are minimal surfaces and are
asymptotic to $R_{\e }^t (\gamma, \sigma , \varsigma)$. Thanks to
Proposition~\ref{quadratique}, we can state that the surface
parameterized by
\[
(1+\gamma )\, (X_\e + u \, N_\e) + \frac{\varsigma}{1+\gamma} \, e_1
+ ( (1+\gamma) \, \log (1+\gamma) + \sigma ) \, e_3
\]
is minimal, if and only if the function $u$ is a solution of
\begin{equation}
L_\e \, u = (1+\gamma)^{-1} \cosh^2 \o  \, Q_\e \left( (1+\gamma )
\frac{u}{\cosh \o } , (1+\gamma ) \, \frac{\nabla u}{\cosh \o},
(1+\gamma ) \, \frac{\nabla^2 u}{\cosh \o}\right) \label{eq:nled}
\end{equation}
The fact that the perturbed surface is asymptotic to the non
perturbed one can then be translated into the fact that the function
$u$ tends to $0$ at $\infty$.

\medskip

We set
\[
\tilde t_\e : = - \frac{1}{2} \, \log \e
\]
Two modifications are now required. First, even though the surface
$R_{\e}^t (\gamma , \sigma , \varsigma)$ can be parameterized by
(\ref{eq:cbos}), its boundary does not correspond to the curve $t =
\tilde t_\e$. We therefore modify the above parametrization so that
the part of $R_{\e}^t ( \gamma , \sigma , \varsigma)$ corresponding
to $t \geq \tilde t_\e + \log 4$ is still parameterized by
(\ref{eq:cbos}), while over the graph over the annulus $B_{2 \,
\e^{-1/2}} - B_{\e^{-1/2}/2}$ for the function defined in
(\ref{eq:vgtt}), we change coordinates
\[
(r, \theta ) = \left( \frac{1}{2} \,  e^t , \tau_\e \, y \right)
\]
Finally, we interpolate (smoothly) between the two parameterizations
in the graph over the annulus $B_{3 \, \e^{-1/2}} - B_{\e^{-1/2}}$
by the function (\ref{eq:vgtt}).

\medskip

The next modification we need to do is concerned with the normal
vector field about $R_{\e}^t  (\gamma , \sigma , \varsigma)$ since
we would like this vector field to be vertical near the boundary of
this surface. This can be achieved by modifying the normal vector
field into a transverse vector field $\tilde N_\e$ which agrees with
the unit normal vector field $N_\e$ for all $t \geq \tilde t_\e +
\log 4$ and which agrees with $e_3$ for all $t \in [\tilde t_\e ,
\tilde t_\e + \log 2]$.

\medskip

Now, we consider a graph over this surface for some function $u$,
using the modified vector field $\tilde N_\e$. This graph will be
minimal if and only if the function $u$ is a solution of some
nonlinear elliptic equation which is not exactly equal to
(\ref{eq:nled}) because of the above two modifications. Indeed,
starting from (\ref{eq:nled}) and taking into account the effects of
the change of parametrization and the change in the vector field
$N_\e$ into $\tilde N_\e$, we see that the minimal surface equation
now reads
\begin{equation}
L_\e  \, u   =  \tilde L_\e \, u + \cdh   \,  \tilde Q_\e \left(
\frac{u}{\cosh \o } , \frac{\nabla u}{\cosh \o }, \frac{\nabla^2
u}{\cosh \o }\right) \label{azaz}
\end{equation}
The nonlinear operator $\tilde Q_\e$ enjoys the same properties as
$Q_\e$ in Proposition~\ref{quadratique}. We will write for short
\[
\hat Q_\e (u) : = \tilde Q_\e \left( \frac{u}{\cosh \o } ,
\frac{\nabla u}{\cosh \o }, \frac{\nabla^2 u}{\cosh \o }\right).
\]
Observe that $\tilde Q_\e$ is explicitly given by
\[
\tilde Q_\e  (\cdot) =  (1+\gamma)^{-1} \, Q_\e ( (1+\gamma) \,
\cdot \, )
\]
when $t \geq \tilde t_\e + \log 4$.

\medskip

The operator $\tilde L_\e$ is a linear second order operator whose
coefficients are supported in $[\tilde t_\e , \tilde t_\e + \log 4]
\times S^1 (\tau_\e)$ and are bounded by  a constant times
$\e^{1/2}$, in ${\mathcal C}^\infty$ topology, where partial
derivatives are computed with respect to the vector fields
$\partial_t$ and $\partial_y$. Let us briefly comment on the
estimate of the coefficients of $\tilde L_\e$. If we were only
taking into account the effect of the change from $N_\e$ into
$\tilde N_\e$, we would obtain, applying the result of Appendix B, a
similar formula where the coefficients of the corresponding operator
$\tilde L_\e$ are bounded by a constant times $\e$  since
\[
\tilde N_\e \cdot N_\e  = 1 + {\mathcal O} (\e)
\]
when $t \in [\tilde t_\e , \tilde t_\e + \log 2]$. If we were only
taking into account the effect of the change in the parametrization,
we would obtain a similar formula where the coefficients of the
corresponding operator $\tilde L_\e$ are bounded by a constant times
$\e^{1/2}$, this basically follows from (\ref{eq:ttop}) which shows
that
\[
U_t^0 (r, \theta ) =  \log (2r) + {\mathcal O} (\e^{1/2}) ,
\]
where the change of coordinates takes place. The estimate of the
coefficients of $\tilde L_\e$ follows from these considerations.

\medskip

Now, assume that we are given a function $\varphi \in {\mathcal
C}^{2, \alpha} (S^1)$ which is even with respect to $y$,
$L^2$-orthogonal to $1$ and $y \longrightarrow \cos (y)$ and which
satisfies
\[
\| \varphi \|_{{\mathcal C}^{2,\a}} \leq \kappa \, \e.
\]
We set
\[
w_\varphi  (\cdot, \cdot) : = {\mathcal P} (\varphi) ( \cdot -
\tilde t_\e , \cdot / \tau_\e )
\]
In order to solve (\ref{azaz}), we choose
\[
\mu \in (-2, -1)
\]
and
look for $u$ of the form
\[
u = w_\varphi +v
\]
where $v \in {\mathcal C}^{2, \alpha}_\mu ([\tilde t_\e , \infty)
\times S^1 (\tau_\e))$. Using the result of
Proposition~\ref{inverse}, we can rephrase this problem as a fixed
point problem
\begin{equation}
v =  S (v) \label{ffp}
\end{equation}
where the nonlinear mapping $S (= S_{\e, \gamma, T_1, \varphi})$
(which depends on $\e, \gamma, T_1$ and $\varphi$) is defined by
\[
S (v) : = {G}_{\e, \tilde t_\e} \left( \tilde L_\e (w_\varphi +v) -
L_\e \, w_{\varphi}  +\cosh ^2 \o \, \hat Q _\e \left( w_\varphi + v
\right) \right).
\]
where the operator $G_{\e, \tilde t_\e}$ is the one defined in
Proposition~\ref{inverse}. The existence of a fixed point of
(\ref{ffp}) is an easy consequence of the following technical~:
\begin{lemma}
\label{dirrie} There exist constants $c_\kappa >0$ and $\e_\kappa >
0$, such that
\begin{equation}
\|S (0) \|_{{\mathcal C}^{2,\a}_\mu} \leq c_\kappa  \, ( \e^{(3
+\mu)/2} + \e^{4+2\mu} ) \label{eq:5.1}
\end{equation}
and, for all $\e \in (0, \e_\kappa)$
\[
\|S (v_2)-S (v_1)\|_{{\mathcal C}^{2,\a}_\mu} \leq \frac{1}{2} \,
\|v_2-v_1 \|_{{\mathcal C}^{2,\a}_\mu}
\]
for all $v_1,v_2 \in {\mathcal C}^{2,\a}_{\mu}([\tilde t_\e ,
\infty) \times S^1 (\tau_\e))$ such that $\| v_i \|_{{\mathcal
C}^{2,\a}_\mu} \leq 2 \, c_\kappa  (\e^{(3+\mu)/2} +  \e^{4+2\mu}
)$.
\end{lemma}
{\bf Proof~:} Using the properties of $w_\varphi$ given in
Lemma~\ref{poissonbis} together with the properties of $\tilde
L_\e$, we immediately get
\[
\| \tilde L_\e (w_\varphi ) \| _{{\mathcal C}^{0,\a}_\mu} \leq
c_\kappa \, \e^{(3+\mu)/2}
\]

Next, we use the fact that
\[
L_\e \, w_\varphi = 2(a^2+b^2) \, w_\varphi
\]
However, we have proved in (\ref{eq:GF}) that  $a^2 \leq \e^2$.
Furthermore, $b^2$ is an even function and, thanks to
Remark~\ref{re:remm}, we know that $b^2 \leq c \, (\cosh t)^{-2}$
for all $t \in [0, t_\e]$ and some constant $c >0$ independent of
$\e$ small enough. Therefore, we conclude (with little work) that
\[
\|L_\e \, w_\varphi \| _{{\mathcal C}^{0,\a}_\mu} \leq  c_\kappa \,
( \e ^{2+\mu/2} + \e^{4+2\mu} )
\]
Observe that the norm on the left hand side is achieved when $t \sim
2 \, t_\e$.

\medskip

 While the last term is easily estimated by
\[
\| \cdh \, \hat  Q _\e \left( w _\varphi \right) \|_{{\mathcal
C}^{0,\a}_\mu} \leq \, c_\kappa  \, \e^{2+\mu/2}
\]
This completes the proof of the first estimate. The second estimate
follows from similar considerations and is left to the reader.
\hfill $\Box$

\medskip

The previous Lemma shows that, provided $\e$ is chosen small enough,
the nonlinear mapping  $S$ is a contraction mapping from the ball of
radius $2 \, c_\kappa (\e^{(3+\mu)/2} + \e^{4+2\mu} )$ in ${\mathcal
C}^{2,\a}_{\mu} ([\tilde t_\e , \infty) \times S^1 (\tau_\e))$ into
itself. Consequently $S$ has a unique fixed point $v$ in this ball.
This provides a minimal surface $R_{\e}^t (\gamma,  \sigma ,
\varsigma , \varphi)$ which is asymptotic to a half Riemann surface
$R_{\e}^t (\gamma, \sigma , \varsigma)$. Observe that near its
boundary, this surface is a vertical graph over the annulus
$B_{\e^{-1/2}} - B_{\e^{-1/2}/2}$ for some function $U_t$ which can
be expanded as
\[
U_t (r, \theta ) = (1+\gamma) \, \log (2 r ) + \sigma  -
\frac{\e}{2} \, r \, \cos \theta - \frac{\varsigma}{r} \, \cos
\theta + {\mathcal P} (\varphi )( \log (2r) -\tilde t_\e  , \theta )
+ V_t(r, \theta)
\]
in which case the boundary of the surface corresponds to $r =
\frac{1}{2} \, \e^{-1/2}$. Here the function $V_t = V_t(\e, \gamma,
\varsigma, \varphi)$ depends nonlinearly on $\gamma, \varsigma$ and
$\varphi$ and satisfies the following
\[
\| V_t(\e, \gamma, \varsigma, \varphi) \|_{{\mathcal
C}^{2,\alpha}_b} \leq \, c \, \e
\]
and
\begin{equation}
\| V_t(\e, \gamma, \varsigma, \varphi) - V_t( \e, \gamma, \varsigma,
\varphi' ) \|_{{\mathcal C}^{2,\alpha}_b} \leq \, c \, (\e^{1/2} +
\e^{3+3\mu/2} ) \, \| \varphi - \varphi'\|_{{\mathcal C}^{2,
\alpha}} \label{eq:contr1}
\end{equation}
where the constant $c >0$ does not depend on $\e$ or $\kappa$ and
$c_\kappa$ only depends on $\kappa$ but not on $\e$. The space
${\mathcal C}^{2, \alpha}_b$ is the space of ${\mathcal
C}^{2,\alpha}$ function where partial derivatives are taken with
respect to the vector fields $r \, \partial_r$ and
$\partial_\theta$.

\medskip

A similar analysis can be carried over starting from the lower end
of Riemann surface to obtain a minimal surface, which will be
referred to as $R_{\e}^b (\gamma, \sigma , \varsigma , \varphi)$,
which is asymptotic to a half Riemann surface and which, near its
boundary is a vertical graph over the annulus $B_{\e^{-1/2}} -
B_{\e^{-1/2}/2}$ for some function $U_b$ which can be expanded as
\[
U_b(r, \theta ) =  - (1+\gamma ) \, \log (2 r) - \sigma -
\frac{\e}{2} \, r \, \cos \theta + \frac{\varsigma}{r} \, \cos
\theta + {\mathcal P}(\varphi) (\log (2r) - \tilde t_\e , \theta ) +
V_b(r, \theta)
\]
in which case the boundary of the surface corresponds to $r =
\frac{1}{2} \, \e^{-1/2}$. The function $V_b$ satisfies  exactly the
same properties as the function $V_t$. Equivalently one can apply a
rotation of angle $\pi$ about the $x_2$-axis to the surface
$R_{\e}^t (\gamma, \sigma , \varsigma , \bar \varphi )$, where $\bar
\varphi ( \cdot ): = - \varphi (\cdot + \pi)$.

\section{An infinite dimensional family of minimal surfaces which are close to $M_k$}

We perform an analysis close to the one performed in the previous
section, starting this time from the minimal surface $M_k (\xi)$
defined in Section 2, for $\xi$ small enough. Recall that the
surface $M_k (\xi)$ has two ends $E_t(\xi)$ and $E_b(\xi)$ which can
be parameterized as in (\ref{eq:paramt}) and (\ref{eq:paramb}). Also
recall that, according to the result of Lemma~\ref{dl3}, a portion
of these ends can be written as a graph over the $x_3=0$ plane for
functions $U_{\xi,t}$ and $U_{\xi,b}$ which are defined in the
annulus $B_{4 \, \e^{-1/2}} - B_{\e^{-1/2}/4}$.

\medskip

Recall that we have defined
\[
\tilde t_\e  = - \frac{1}{2}\, \log \e
\]
As in the previous section, we modify the parametrization of
the end $E_t(\xi)$ which is given by (\ref{eq:paramt}), say when $s
\in [\tilde t_\e - \log 8, \tilde t_\e + \log 4]$, so that, when $r
\in [\e^{-1/2}/ 4 , 2 \, \e^{-1/2} ]$ the curve corresponding to the
image of
\[
\theta \longrightarrow (r\, \cos \theta, r\, \sin \theta, U_{t, \xi
} (r, \theta))
\]
corresponds to the curve $s = \log (2 r)$. We perform a similar task
for the parametrization of $E_b(\xi)$ so that, when $r \in
[\e^{-1/2}/ 4 , 2 \, \e^{-1/2} ]$ the curve corresponding to the
image of
\[
\theta \longrightarrow (r\, \cos \theta, r\, \sin \theta, U_{b ,
\xi} (r, \theta) )
\]
corresponds to the curve $s = - \log (2 r)$.

\medskip

This being understood, as in the previous section, we modify the
unit normal vector field on $M_k (\xi)$ to produce a transverse unit
vector field $\tilde N_\xi$ which coincides with the normal vector
field $N_\xi$ on $M_k (\xi)$, is equal to $e_3$ on the graph over
$B_{2 \, \e^{-1/2}} - B_{3\e^{-1/2}/8}$ of the functions $U_{t, \xi
}$ and $U_{b, \xi}$ and interpolate smoothly in between the
different definitions of $\tilde N_\e$ in different subsets of
$M_k(\xi)$.

\medskip

A graph of the function $u$, using the vector field $\tilde N_\xi$,
will be a minimal surface if and only if $u$ is a solution of a
second order nonlinear elliptic equation of the form
\[
{\mathbb L}_{M_k(\xi)} \, u = \tilde L_{\xi,\e} \, u + Q_{\xi,\e} \,
(u)
\]
where ${\mathbb L}_{M_k(\xi)}$ is the Jacobi operator about $M_k
(\xi)$, $Q_{\xi,\e}$ is a nonlinear second order differential
operator which collects all the nonlinear terms and $\tilde
L_{\xi,\e}$ is a linear operator which take into account the change
of parametrization and the change of the normal vector field $N_\xi$
into $\tilde N_\xi$, which are described above. Now, we can be more
precise and, at the ends $E_t(\xi)$ and $E_b(\xi)$. For example at
$E_t(\xi)$, and granted the above parametrization, the nonlinear
operator $Q_{\xi,\e}$ can be expanded as
\[
Q_{\xi,\e} (u) =  Q_{2, \xi, \e} \left( \frac{w}{\cosh s}\right) +
\cosh s \, Q_{3, \xi, \e} \left( \frac{w}{\cosh s} \right)
\]
where $Q_{2,\xi , \e}$ and $Q_{3 ,\xi, \e}$ are nonlinear second
order differential operators which satisfy (\ref{eq:ZZZ}), uniformly
in $\xi$ and $\e$.

\medskip

The operator $\tilde L_{\xi,\e}$ is a linear operator which is supported in
$[\tilde t_\e - \log 8, \tilde t_\e + \log 4] \times S^1$ and has
coefficients which are bounded by a constant times $\e^{3/2}$,
uniformly in $\xi$ and $\e$ (The rational for this estimate is that
$\e^{3/2} = \e \, \e^{1/2}$, the first $\e$ comes from the conformal
factor $(\cosh s)^{-2}$ and the $\e^{1/2}$ comes from the
modification in the parametrization and the vector field as in the
previous section).

\medskip

Finally, observe that still in $E_t(\xi)$ the difference
\[
(\cosh s)^{4} \, \left({\mathbb L}_{M_k(\xi)} -  \frac{1}{\cosh^2 s}
\, ( \partial_s^2 + \partial_\theta^2 ) \right)
\]
is a second order differential operator in $\partial_s$ and
$\partial_\theta$ whose coefficients are bounded in ${\mathcal
C}^\infty$ topology (where partial derivatives are taken with
respect to the vector fields $\partial_s$ and $\partial_\theta$)
uniformly in $\xi$ and $\e$. All these facts follow from the
expansion provided in (\ref{eq:2.2}).

\medskip

Now, assume that we are given two functions $\varphi_t, \varphi_b
\in {\mathcal C}^{2, \alpha} (S^1)$ which is even with respect to
$\theta$ and $L^2$ orthogonal to $1$ and $\theta \longrightarrow
\cos \, \theta$ and satisfy
\[
\|\varphi_t\|_{{\mathcal C}^{2, \alpha}} + \|\varphi_b\|_{{\mathcal
C}^{2, \alpha}} \leq \kappa \, \e.
\]
We set $\Phi : = (\varphi_t, \varphi_b)$ and we define $w_\Phi$ to
be the function which is equal to $\chi_t \, {\mathcal P}
(\varphi_t) (\cdot - t_\e , \cdot)$ on the image of $X_{t, \xi}$
where $\chi_t$ is a cutoff function equal to $0$ for $s \leq s_0+1$
and identically equal to $1$ for $s\geq s_0+2$, and is equal to
$\chi_b \, {\mathcal P} (\varphi_b) ( \cdot + t_\e , \cdot)$ on the
image of $X_{b, \xi}$ where $\chi_b$ is a cutoff function equal to
$0$ for $s \geq -s_0-1$ and identically equal to $1$ for $s\leq
-s_0-2$.

\medskip

We define $M_k (\xi, \e)$ to be equal to $M_k (\xi)$ with the image
of $(\tilde t_\e, +\infty) \times S^1$ by $X_{t,\xi}$ and the image
of $(-\infty, - \tilde t_\e) \times S^1$ by $X_{b,\xi}$ removed. We
would like to solve the equation
\[
{\mathbb L}_{M_k(\xi)} (w_\Phi + v)  =\tilde  L_{\xi,\e} (w_\Phi + v) +
Q_{\xi, \e} (w_\Phi + v )
\]
on $M_k(\xi, \e)$, so that the graph of $w_\Phi +v$ will be a
minimal surface.

\medskip

We choose
\[
\delta \in (1,2)
\]
and use the result of Proposition~\ref{pr:lk} so that we can
rephrase the above problem as a fixed point problem
\begin{equation}
v= T(v)
\label{eq:sdsd}
\end{equation}
where
\[
T(v) = {G}_{\xi, \delta} \circ {\mathcal E}_\e \,  \left(
\tilde L_{\xi,\e} (w_\Phi  + v) + {\mathbb L}_{M_k(\xi)} \, w_\Phi +
Q_{\xi,\e} (w_\Phi + v ) \right)
\]
where ${\mathcal E}_\e$ is an extension (linear) operator
\[
{\mathcal E}_\e : {\mathcal C}^{0, \alpha}_{\delta-2, 4}( M_{k}(\xi,
\e)) \longrightarrow {\mathcal C}^{0, \alpha}_{\delta-2,4}
(M_k(\xi)) ,
\]
defined by ${\mathcal E}_\e \, v = v$ in $M_k (\xi, \e)$, ${\mathcal
E}_\e \, v =0$ on the image of $[\tilde t_\e + 1, +\infty) \times
S^1$ by $X_{t,\xi}$ and the image of $(-\infty, - \tilde t_\e -1]
\times S^1$ by $X_{b,\xi}$ and ${\mathcal E}_\e \, v$ interpolate
between these so that, for example,
\[
({\mathcal E}_\e v )\circ X_{t,\xi} (t, \theta )= (1+ \tilde t_\e
-t) \,  v )\circ X_{t,\xi} (\tilde t_\e, \theta )
\]
for $(t, \theta) \in [\tilde t_\e, \tilde t_\e+1]\times S^1$. Here
${\mathcal C}^{0, \alpha}_{\delta-2, 4}( M_{k}(\xi, \e))$ is the
space of restrictions of elements of ${\mathcal C}^{0,
\alpha}_{\delta-2,4} (M_k(\xi))$ to $M_k(\xi, \e)$, endowed with the
induces norm.

\medskip

As in Section 5, the existence of a fixed point $v \in {\mathcal
C}^{2, \alpha}_{\delta, 0} (M_k (\xi))$ for (\ref{eq:sdsd}) follows
at once from the technical~:
\begin{lemma}
\label{dircosta} There exist constants $c_\kappa  >0$ and $\e_\kappa
>0$, such that
\begin{equation}
\|T(0) \|_{{\mathcal C}^{2,\a}_{\delta,0}} \leq c_\kappa \, \e^2
\label{eq:A}
\end{equation}
and, for all $\e \in (0, \e_\kappa)$
\[
\|T (v_2)-T(v_1)\|_{{\mathcal C}^{2,\a}_{\delta,0}} \leq \frac{1}{2}
\, \| v_2-v_1\|_{{\mathcal C}^{2,\a}_{\delta,0}}
\]
for all $v_1,v_2 \in {\mathcal C}^{2,\a}_{\delta,0} (M_k(\xi))$ and
satisfying $\|v\|_{{\mathcal C}^{2,\a}_{\delta,0}} \leq \, 2 \,
c_\kappa \, \e^2$.
\end{lemma}
\begin{proof}
The proof is similar to the one in the proof of Lemma~\ref{dirrie}.
Again, we use the result of Lemma~\ref{poissonbis} to obtain the
estimate
\[
\| {\mathcal E}_\e ( {\mathbb L}_{M_k(\xi)} \, w_\Phi )
\|_{{\mathcal C}^{0,\a}_{\delta -2, 4}} \leq c_\kappa \, \e^2
\]
and, using the properties of $\tilde L_{\xi , \e}$, we obtain
\[
\| {\mathcal E}_\e (\tilde L_{\xi,\e} \, w_\Phi )\|_{{\mathcal
C}^{0,\a}_{\delta -2, 4}} \leq c_\kappa \, \e^{\frac{3+\delta}{2}}
\]
Finally, we have
\[
\| {\mathcal E}_\e (Q_{\xi,\e}  (w_\Phi)) \|_{{\mathcal C}^{0,
\a}_{\delta-2, 4}} \leq  c_\kappa \, \e^{3+\delta/4}
\]
We leave the details to the reader. \end{proof}

\noindent

The previous Lemma shows that, provided $\e$ is chosen small enough,
the nonlinear mapping $T$ is a contraction mapping from the ball of
radius $2 \, c_\kappa \,\e^2$ in ${\mathcal C}^{2,\a}_{\delta, 0}
(M_k (\xi))$ into itself. Consequently $T$ has a unique fixed point
$v$ in this ball. This provides a minimal surface $M_k (\xi, \e,
\varphi_t, \varphi_b)$ which is close to $M_k (\xi, \e)$, has one
horizontal end and two boundaries. This surface is, close to its
upper boundary, a vertical graph over the annulus $B_{\e^{-1/2}/2} -
B_{\e^{-1/2}/4}$ for some function $\bar U_t$ which can be expanded
as
\[
\bar U_t (r, \theta ) = \sigma_{t, \xi} + \log (2 r) + \xi \, r \,
\cos \theta + {\mathcal P} (\varphi_t) (\tilde t_\e - \log (2r) ,
\theta ) + \bar V_t(r, \theta)
\]
and this surface is, close to its lower boundary, a vertical graph
over the annulus $B_{\e^{-1/2}/2} - B_{\e^{-1/2}/4}$ for some
function $\bar U_b$ which can be expanded as
\[
\bar U_b (r, \theta ) = - \sigma_{b, \xi} - \log (2 r) + \xi \, r \,
\cos \theta + {\mathcal P} (\varphi_b) (\tilde t_\e - \log (2r) ,
\theta ) + \bar V_b (r, \theta)
\]
where $\bar V_ t =  \bar V_t (\e, \xi, \Phi)$ and $\bar V_ b =  \bar
V_b (\e, \xi, \Phi)$ depend nonlinearly on $\e , \xi$ and $\Phi$ and
satisfy (for $V=V_t$ or $V=V_b$)
\[
\| \bar V (\e, \xi , \Phi) \|_{{\mathcal C}^{2,\alpha}_b} \leq \, c
\, \e
\]
and
\begin{equation}
\| \bar V (\e, \xi  , \Phi) - \bar V ( \e, \xi , \Phi' )
\|_{{\mathcal C}^{2,\alpha}_b} \leq \, c \, \e^{1 - \delta/2} \, \|
\Phi - \Phi'\|_{{\mathcal C}^{2, \alpha}} \label{eq:contr2}
\end{equation}
where the constant $c
>0$ does not depend on $\e$ or $\kappa$ and $c_\kappa$ only depends
on $\kappa$ but not on $\e$. The boundaries of the surface
corresponds to $r = \frac{1}{2} \, \e^{-1/2}$.

\section{The matching of Cauchy data and the proof of the main result}

We collect the results we have obtained in Section 5 and Section 6.
In Section 5, we have obtained two surfaces which are perturbations
of the upper (rep. the lower end) of Riemann's surface. The first
surface
\[
R_{\e+\eta_t}^t (\gamma_t, \sigma_{t, \xi} + \sigma_{t} ,
\varsigma_{t} , \varphi_t )
\]
depends on the parameters $\eta_t, \gamma_t, \sigma_{t},
\varsigma_{t}$ and the function $\varphi_t$ and can be
parameterized, close to its boundary as the vertical graph of
\[
U_t(r, \theta) : = (1+ \gamma_t) \, \log ( 2r ) + \sigma_{t,\xi} +
\sigma_t - \frac{\e+\eta_t}{2} \, r \, \cos \theta -
\frac{\varsigma_t}{r} \, \cos \theta + {\mathcal P}(\varphi_t) (\log
r - \tilde t_\e , \theta)  + {\mathcal O} (\e)
\]
The second surface
\[
R_{\e+\eta_b}^b (\gamma_b, \sigma_{b, \xi} + \sigma_b, \varsigma_{b}
, \varphi_b)
\]
depends on the parameters $\eta_b, \gamma_b, \sigma_b, \varsigma_b$
and the function $\varphi_b$ and can be parameterized, close to its
boundary as the vertical graph of
\[
U_b(r, \theta) : =  - (1+ \gamma_b) \, \log  (2r) - \sigma_{b,\xi} -
\sigma_b - \frac{\e + \eta_b}{2} \, r \, \cos \theta +
\frac{\varsigma_b}{r} \, \cos \theta + {\mathcal P} (\varphi_b)
(\log r - \tilde t_\e, \theta)  + {\mathcal O} (\e)
\]

Now, collecting the result of Section 6, we have a surface
\[
M_k (\xi, \e, \tilde \varphi_b, \tilde \varphi_t)
\]
which has two boundaries, one end asymptotic to a horizontal plane
and can be parameterized, close to its upper boundary as the
vertical graph of
\[
\bar U_t(r, \theta) : =  \log (2r)  + \sigma_{t,\xi} + \xi \, r\,
\cos \theta + {\mathcal P}(\tilde \varphi_t) (\tilde t_\e- \log r,
\theta) + {\mathcal O} (\e)
\]
while it can be parameterized close to its lower boundary as the
vertical graph of
\[
\bar U_b(r, \theta) :=  - \log (2r) - \sigma_{b,\xi} + \xi \, r\,
\cos \theta + {\mathcal P}(\tilde \varphi_b) (\tilde t_\e - \log r,
\theta) + {\mathcal O} (\e)
\]

We set \[
\xi = - \frac{\e}{2}
\]
and assume that the parameters and the boundary functions are chosen
so that
\[
\begin{array}{lllll}
\e^{-1/2} \, (|\eta_t| + |\eta_b| ) + \e^{1/2} \, (|\varsigma_t|+
|\varsigma_b|) +  |\log \e|^{-1} \,(|\sigma_t|+ |\sigma_b|) \\[3mm]
\qquad \qquad +  |\gamma_{t}|+ |\gamma_t| +  \|
\varphi_t\|_{{\mathcal C}^{2, \alpha}} + \| \varphi_b\|_{{\mathcal
C}^{2, \alpha}}+ \| \tilde \varphi_t\|_{{\mathcal C}^{2, \alpha}}+
\| \tilde \varphi_b\|_{{\mathcal C}^{2, \alpha}} \leq \kappa \, \e
\end{array}
\]
where the constant $\kappa >0$ is fixed large enough. Recall that
the functions $\varphi_t, \varphi_b, \tilde \varphi_t$ and $\tilde
\varphi_b$ are assumed to be even and $L^2$ orthogonal to the
functions $1$ and $\theta \longrightarrow \cos \theta$. The
functions ${\mathcal O}(\e)$ do depend nonlinearly on the different
parameters and boundary data functions but are bounded by a constant
(independent of $\kappa$ and $\e$) times $\e$ in ${\mathcal C}^{2,
\alpha}$ topology, when partial derivatives are taken with respect
to the vector fields $r \, \partial_r$ and $\partial_\theta$.

\medskip

It remains to show that, for all $\e$ small enough,  it is possible
to choose the parameters and boundary functions in such a way that
the surface
\[
R_{\e+\eta_t}^t (\gamma_t, \sigma_{t, \xi} + \sigma_{t} ,
\varsigma_{t} , \varphi_t ) \cup M_k (\xi, \e, \tilde \varphi_b,
\tilde \varphi_t) \cup R_{\e+\eta_b}^b (\gamma_b, \sigma_{b, \xi} +
\sigma_b, \varsigma_{b} , \varphi_b)
\]
is a ${\mathcal C}^1$ surface across the boundaries of the different
summands. Regularity theory will then ensure that this surface is in
fact smooth and by construction is has the desired properties. This
will therefore complete the proof of the main theorem.

\medskip

Granted the description of the surfaces close to their respective
boundaries it is enough to fulfill that following system of
equations
\[
U_t = \bar U_t \qquad U_b = \bar U_b \qquad \partial_r U_t =
\partial_r \bar U_t \qquad \partial_r U_b =
\partial_r \bar U_b
\]
on $S^1( \frac{1}{2} \, \e^{-1/2})$.

\medskip

The first two equations lead to the system
\begin{equation}
\begin{array}{rllll}
- \frac{1}{2} \, \log \e \, \gamma_t + \sigma_t - \left( \frac{1}{4}
\, \e^{-1/2} \, \eta_t  + 2 \, \e^{1/2} \, \varsigma_{t} \right) \,
\cos \theta  + \varphi_t  - \tilde \varphi_t& = & {\mathcal O}(\e)\\[3mm]

\frac{1}{2} \, \log \e \, \gamma_b - \sigma_{b} - \left( \frac{1}{4}
\, \e^{-1/2} \, \eta_b  - 2 \, \e^{1/2} \, \varsigma_{b} \right) \,
\cos \theta + \varphi_b - \tilde \varphi_b & = & {\mathcal O}(\e)
\end{array}
\label{eq:dd1}
\end{equation}
while the last two equations read
\begin{equation}
\begin{array}{rllll}
\gamma_t  - \left( \frac{1}{4} \, \e^{-1/2} \, \eta_t - 2 \,
\e^{1/2} \, \varsigma_{t} \right) \, \cos
\theta + \partial_t ({\mathcal P} (\varphi_t+ \tilde \varphi_t)) & = & {\mathcal O}(\e)\\[3mm]

- \gamma_b  - \left( \frac{1}{4} \, \e^{-1/2} \, \eta_b + 2 \,
\e^{1/2} \, \varsigma_{b} \right) \, \cos \theta  + \partial_t
({\mathcal P} (\varphi_b+\tilde \varphi_b))& = & {\mathcal O}(\e)
\end{array}
\label{eq:dd2}
\end{equation}

Projection of every equation of this system over the
$L^2$-orthogonal complement of $\mbox{Span} \{1, \cos\}$, we obtain
the system
\[
\begin{array}{rlllrlrrr}
 \varphi_t  - \tilde \varphi_t& = & {\mathcal O}(\e) & \qquad & \varphi_b - \tilde \varphi_b & = & {\mathcal O}(\e) \\[3mm]
\partial_t {\mathcal P}(\varphi_t  +  \varphi_t) & = & {\mathcal O}(\e) & \qquad & \partial_t {\mathcal P}(\varphi_b+ \varphi_b)& = & {\mathcal O}(\e)
\end{array}
\]
Observe that the operator
\[
\begin{array}{clclll}
{\mathcal C}^{2, \alpha} & \longrightarrow & {\mathcal C}^{1,
\alpha}\\[3mm]
\varphi & \longmapsto & \partial_t {\mathcal P}(\varphi)
\end{array}
\]
is invertible, and hence the last system can be rewritten as
\begin{equation}
\left(  \varphi_t, \tilde \varphi_t, \varphi_b , \tilde \varphi_b
\right) =  {\mathcal O}(\e) \label{eq:qq}
\end{equation}
Recall that the right hand side depends nonlinearly on $ \varphi_t,
\tilde \varphi_t, \varphi_b , \tilde \varphi_b$. Thanks to
(\ref{eq:contr1}) and (\ref{eq:contr2}) we can use a fixed point
theorem for contraction mapping in the ball of radius $\kappa \, \e$
in $({\mathcal C}^{2, \alpha} )^4$ to obtain, for all $\e$ small
enough, a solution (\ref{eq:qq}) which depends at least continuously
(and in fact smoothly) on the parameters $\gamma_t, \gamma_b ,
\sigma_{t}, \sigma_{b} , \varsigma_{t} , \varsigma_{b} , \eta_t $
and $\eta_b$.

\medskip

Inserting this solution into (\ref{eq:dd1}) and (\ref{eq:dd2}), we
see that it remains to solve a system of the form
\begin{equation}
\begin{array}{rllll}
- \frac{1}{2} \, \log \e \, \gamma_t + \sigma_t - \left( \frac{1}{4}
\, \e^{-1/2} \, \eta_t  + 2 \, \e^{1/2} \, \varsigma_{t} \right) \,
\cos \theta  & = & {\mathcal O}(\e)\\[3mm]
\frac{1}{2} \, \log \e \, \gamma_b - \sigma_{b} - \left( \frac{1}{4}
\, \e^{-1/2} \, \eta_b  - 2 \, \e^{1/2} \, \varsigma_{b} \right) \,
\cos \theta  & = & {\mathcal O}(\e) \\[3mm]
\gamma_t  - \left( \frac{1}{4}
\, \e^{-1/2} \, \eta_t - 2 \, \e^{1/2} \, \varsigma_{t} \right) \,
\cos \theta  & = & {\mathcal O}(\e)\\[3mm]
- \gamma_b  - \left( \frac{1}{4} \, \e^{-1/2} \, \eta_b + 2 \,
\e^{1/2} \, \varsigma_{b} \right) \, \cos \theta & = & {\mathcal
O}(\e)
\end{array}
\end{equation}
where the right hand sides depend nonlinearly on $\gamma_t, \gamma_b
, \sigma_{t}, \sigma_{b} , \varsigma_{t} , \varsigma_{b} , \eta_t $
and $\eta_b$.

\medskip

Projecting this system over the constant function and the function
$\theta \longrightarrow \cos \theta$, we see that this system can be
rewritten as
\begin{equation}
\left( \gamma_t, \gamma_b ,   \bar \sigma_{t}, \bar \sigma_{b} ,
\bar \varsigma_{t} , \bar \varsigma_{b} ,  \bar \eta_t , \bar \eta_b
\right) = {\mathcal O}(\e) \label{eq:qqq}
\end{equation}
where we have set
\[
(\bar \varsigma_{t}, \bar \varsigma_{b}) := \e^{1/2} \,
(\varsigma_{t}, \varsigma_{b}), \qquad  (\bar \sigma_{t}, \bar
\sigma_{b}) : = |\log \e|^{-1} \, (\sigma_{t} , \sigma_{b})
\]
and
\[
(\bar \eta_t ,  \bar \eta_b) : = \e^{-1/2} \, (\eta_t , \eta_b)
\]
This time we can use Leray-Schauder degree theory in the ball of
radius $\e$ in $\R ^8 $ to solve (\ref{eq:qqq}), for all $\e$ small
enough. This completes the proof of a solution of
(\ref{eq:dd1})-(\ref{eq:dd2}) and hence the proof of the main
theorem.

\medskip

\begin{remark}
Alternatively, with more work, one can use a fixed point argument
for contraction mapping  to solve (\ref{eq:qqq}).
\end{remark}

\section{Appendix A}

We consider the surface parameterized by
\[
X  = X_c + w \, N_c
\]
The coefficients of $g_w$, the first fundamental form of this
surface, are given by
\[
|\partial_s X|^2 =  \cosh^2 s - 2 \, w + \frac{1}{\cosh^2 s} \, w^2
+ (\partial_s w)^2
\]
\[
|\partial_\theta  X|^2 =  \cosh^2 s + 2 \, w + \frac{1}{\cosh^2 s} \, w^2
+ (\partial_\theta w)^2
\]
and
\[
\partial_s X \cdot \partial_\theta X =  \partial_s w \, \partial_\theta w
\]
It follows from these that the determinant of the metric $g_w$ can
be expanded as
\[
\begin{array}{rlllll}
|g_w| & = & \displaystyle \cosh^4 s \, \left( 1 + \frac{1}{\cosh^2
s} \, \left( (\partial_s w)^2 + (\partial_\theta w)^2 -
\frac{2}{\cosh^2 s} \, w^2 \right) \right. \\[3mm]
& + & \displaystyle \left( \frac{1}{\cosh s} \, P_3 (\frac{w}{\cosh
s} , \frac{\nabla w}{\cosh s}) + P_4 (\frac{w}{\cosh s} ,
\frac{\nabla w}{\cosh s}) \right)
\end{array}
\]
where the $P_i$ are homogeneous polynomials of degree $i$, whose
coefficients are bounded smooth functions of $s$ and $\theta$.

\medskip

We consider the area energy
\[
A (w) := \int \sqrt{|g_w|} \, ds\, d\theta
\]
The surface parameterized by $X_w$ is minimal if and only if the
first variation of $A$ at $w$, is $0$. This can be written as
\[
2 \, DA_{|w} (v) = \int \frac{1}{\sqrt{|g_w|}} \, D_w |g_w| \, (v)
\, ds \, d\theta
\]
Observe that
\begin{equation}
\begin{array}{rllll} \frac{1}{\sqrt{|g_w|}} \, D |g_w|_{_w}
\, (v) & = & \displaystyle  \partial_s w \, \partial_s v +
\partial_\theta w \, \partial_\theta v - \frac{2}{\cosh^2 s}
\, w \, v  \\[3mm]
&+& \displaystyle \left( \tilde Q_2 \left( \frac{w}{\cosh s },
\frac{ \nabla w}{\cosh s } \right) + \cosh s \, \tilde Q_3 \left(
\frac{w}{\cosh s }, \frac{ \nabla w}{\cosh s} \right)  \, \right) \,
v \\[3mm]
&+& \displaystyle  \left( \tilde Q_2' \left( \frac{w}{\cosh s },
\frac{ \nabla w}{\cosh s } \right) + \cosh s \, \tilde Q'_3 \left(
\frac{w}{\cosh s }, \frac{ \nabla w}{\cosh s} \right)  \, \right) \,
\partial_s v \\[3mm]
&+& \displaystyle  \left( \tilde Q_2'' \left( \frac{w}{\cosh s },
\frac{ \nabla w}{\cosh s } \right) + \cosh s \, \tilde Q''_3 \left(
\frac{w}{\cosh s }, \frac{ \nabla w}{\cosh s} \right)  \, \right) \,
\partial_\theta v
\end{array}
\label{eq:pippppplll}
\end{equation}
where the operator $Q_2, \ldots $ and the operators $\tilde Q_3,
\ldots $ enjoy properties similar to the one enjoyed by $Q_2$ and
$Q_3$ in the statement of the result. The result then follows at
once.

\section{Appendix B}

This appendix is essentially a generalization of the corresponding
analysis in \cite{Maz-Pac}. Let be $\Sigma$ be a smooth surface
embedded in a Riemannian manifold $(M, g)$. We denote by $N$ the
unit normal vector field compatible with the orientation of
$\Sigma$. Suppose that $\tilde N$ is another unit vector field
transverse to $\Sigma$, the implicit function theorem implies that,
given $p_0 \in \Sigma$, there exist neighborhoods ${\mathcal U}$ and
${\mathcal V}$ of $(p_0 ,0) \in \Sigma \times {\mathbb R}$ and a
diffeomorphism $(p,s) \longrightarrow (\varphi (p,s),\psi (p,s))$
from ${\mathcal U}$ to ${\mathcal V}$ such that
\begin{equation}
\mbox{Exp}^M_p (s \, \tilde N (p)) = \mbox{Exp}^M_{\varphi (p,s)}
(\psi (p,s) \, N (\varphi (p,s) )) \label{eq:a1}
\end{equation}
where $\mbox{Exp}^M$ denotes the exponential map in $(M,g)$. In
addition $\varphi (p,0)=p$ and $\psi (p,0)=0$.

\medskip

Differentiation of (\ref{eq:a1}) with respect to $s$ at $s=0$ yields
\begin{equation}
\tilde N (p) = \partial _s \varphi (p,0) + \partial _s \psi (p,0) \,
N(p)
\label{pipo}
\end{equation}
Taking the scalar product with $N(p)$ we conclude that
\begin{equation}
g( \tilde N (p), N(p) ) = \partial _s \psi (p,0) \label{eq:pdpsi}
\end{equation}
This immediately implies that $\psi (p,s) = g ( \tilde N (p), N(p))
\, s + {\mathcal O}(s^2)$. On the other hand, projection of
(\ref{pipo}) over $T_p\Sigma$ yields
\begin{equation}
\tilde N^t (p) = \partial _s \varphi (p,0) \label{eq:pdphi}
\end{equation}
where $\tilde N^t (p)$ is the tangential component of $\tilde N$.

\medskip

Next any surface $\tilde \Sigma$ sufficiently close to $\Sigma$ can
be either parameterized as a graph of the function $w$ over $\Sigma$
using the vector field $\tilde N$ or the graph of the function $\bar
w$ for the normal vector field $N$. Thanks to the above analysis we
can write
\[
\bar w (\varphi (p, w(p))) =  \psi (p, w(p))
\]
Now, the mean curvature of the surface $\Sigma$ at the point
$\mbox{Exp}^M_p (w(p) \, \tilde N (p))$ and at the point
$\mbox{Exp}^M_{\bar p} ( \bar w (\bar p) \, N (\bar p ))$ are the
same if $ \bar p =  \varphi (p, w(p))$. We phrase this property as
\[
H_{\tilde N , w} (p) =  H_{N, \bar w} (p)
\]
Differentiation with respect to $w$ at $w=0$ yields
\[
DH_{\tilde N , 0}  =  DH_{N, 0} (\partial_s \psi ) +
\nabla_{\partial_s \varphi}  H_{N,0}
\]
Taking into account the partial derivatives of $\varphi$ and $\psi$,
which are given in (\ref{eq:pdphi}) and (\ref{eq:pdpsi}), we
conclude that
\[
DH_{\tilde N , 0} \, u =  DH_{N, 0} ( g( \tilde N (p), N(p) )  \, u
) + \left(\nabla_{\tilde N^t (p)} H_{N,0} \right)\, u
\]
for any smooth function $u$ defined on $\Sigma$. In the special case
where $\Sigma$ has constant mean curvature, we simply get
\[
DH_{\tilde N , 0} \, (u) =  {\mathcal L}_{\Sigma} ( g( \tilde N (p),
N(p) ) \, u )
\]
which gives the relation between ${\mathcal L}_{\Sigma}$ the Jacobi
operator about $\Sigma$ and $DH_{\tilde N , 0}$ the linearized mean
curvature operator when the normal vector field $N$ is changed into
a transverse vector field $\tilde N$.

\end{document}